\documentclass[a4paper,11pt]{article}

\usepackage[centertags]{amsmath}
\usepackage[nosort]{cite}
\usepackage[psamsfonts]{amsfonts} 
\usepackage{amsthm,enumerate,amssymb}

\newcommand{\dd}{\mathrm{d}}
\newcommand{\im}{\mathrm{i}}
\newcommand{\e}{\mathrm{e}}
\newcommand{\bigtimes}{\mathop{\mbox{\Large $\times$}}} 
\newcommand{\QED}{\hfill$\square$ \\ \medskip}

\newcommand{\bbbone}{{\mathchoice {\rm 1\mskip-4mu l} {\rm 1\mskip-4mu l}
{\rm 1\mskip-4.5mu l} {\rm 1\mskip-5mu l}}}

\newcommand{\ess}{\Gamma_\mathbb{R}^G}
\newcommand{\essp}{\Gamma_{\mathbb{R}_+}^G}
\newcommand{\essn}{\Gamma_{(-\infty,0]}^G}
\newcommand{\is}{\mathcal{H}}
\newcommand{\ms}{\mathfrak{Z}}

\theoremstyle{plain}
\newtheorem{theorem}{Theorem}

\newtheorem{proposition}[theorem]{Proposition}

\begin{document}

\title{Classical dilations \`a la Hudson-Parthasarathy of Markov semigroups}

\author{M.~Gregoratti\thanks{E-mail:
matteo.gregoratti@polimi.it}\\ \\
{\footnotesize \textsl{Dipartimento di Matematica ``F.Brioschi'', Politecnico di
Milano}} \\
{\footnotesize \textsl{Piazza Leonardo da Vinci 32, I-20133 Milano, Italy}}\\ \\
{\small Quaderno di Dipartimento \textbf{QDD 15}}}

\date{}

\maketitle

\begin{abstract}
\noindent We study the Classical Probability analogue of the dilations of 
a quantum dynamical semigroup defined in Quantum Probability via quantum stochastic differential equations. 
Given a homogeneous Markov chain in continuous time in a finite state space $E$,  
we introduce a second system, an environment, and a 
deterministic invertible time-homogeneous global evolution of the system $E$ with this environment
such that the original Markov evolution of $E$ can be realized by a proper choice of the initial random 
state of the environment. 
We also compare this dilations with the dilations of a quantum dynamical semigroup in Quantum Probability: given a classical Markov semigroup,
we extend it to a proper quantum 
dynamical semigroup for which we can find a Hudson-Parthasarathy dilation which is itself an extension of our classical dilation.

\bigskip

\noindent AMS Subject Classification: 60J27, 81S25
\end{abstract}

\section{Introduction}

We study the analogue in Classical Probability of the dilations in Quantum Probability of a quantum dynamical 
semigroup (QDS) in continuous time. 
A QDS $T_t$ describes the evolution of a quantum system, possibly open, but
``Markovian'', and homogeneous in time. If a QDS is uniformly continuous, then it is always
possible to introduce a
Hudson-Parthasarthy quantum stochastic differential equation and to employ its solution to dilate $T_t$ by a
quantum stochastic flow $j_t$. Such a dilation allows to represent the QDS by the conditional expectation of 
a quantum Markov process, analogously to
the representation of a classical Markov semigroup (CMS) by a classical Markov process.
Anyway, such a dilation in Quantum
Probability enjoys a reacher structure which allows to dilate the semigroup $T_t$ at the same time
also by a strongly continuous unitary group $U_t$, thus showing that the system evolution $T_t$ and the flow
$j_t$ do not contradict the axioms of Quantum Mechanics, i.e.\ that they can arise from a Hamiltonian 
evolution $U_t$ of an isolated bigger system, consisting of the given system and its 
environment \cite{ALV,F1,F2,HP,M1,M2,M,MS,P}. 
In particular this implies that the Hamiltonian operator generating $U_t$ gives
an infinitesimal description of the dilation which is alternative, but equivalent, to that given by the
Hudson-Parthasarathy equation \cite{C1,C2,bm1,bm2}. This feature differentiates the representation of a QDS
by means of a quantum stochastic flow from the representation of a CMS by means of a Markov process, as
usually there is no ``Hamiltonian evolution'' associated to this latter.

Dilations analogous to the quantum ones have been recently introduced in Classical Probability for CMS in
discrete time \cite{bm7}. The first aim of this paper is to introduce them also in the continuous time context, 
choosing a
self-contained approach in a completely classical framework. These classical dilations are interesting, not
only to better understand the relationship between the two probabilistic theories, but also from a simply
classical point of view, to better understand the relationship between Markov processes and 
deterministic invertible homogeneous dynamics.
Indeed, we show that every homogeneous Markov chain in continuous time in a finite state space can
always be realized as a deterministic invertible homogeneous evolution 
of the system coupled with a second system.
The existence of such representations is theoretically relevant if Markov chains
are applied to phenomena, like physical phenomena for example, for which an underlaying theory
postulates deterministic invertible homogeneous evolutions in absence of noise and external
disturbances. For these phenomena the second system introduced by the dilation models the surrounding
world, the environment, the source of the noise, which is given now a dynamical explanation. 
Of course, here the characterizing property is not simply that the global evolution is deterministic, as in
innovation theory \cite{R,T}, but that it is also invertible and homogeneous in time.

More precisely, we consider a system with finite state space $E$, undergoing a continuous 
time evolution given by a homogeneous Markov chain. 
Then, we introduce an environment  
with its state space $(\Gamma, \mathcal{G})$, a measurable space, together with a group of
measurable maps
$\alpha_t:E\times\Gamma\to E\times\Gamma$, $t\in\mathbb{R}$, describing a global evolution.
Thus, if $(i,\gamma)$ is the state of the
compound system at time 0, then $\alpha_t(i,\gamma)$ is its state at time $t$,
where hence $\alpha_t$ gives a deterministic invertible homogeneous global evolution. Nevertheless, 
if the environment state is never observed and if initially it is randomly
distributed with some law $\mathbb{Q}$ on $(\Gamma,\mathcal{G})$, then the evolution of the observed system turns 
out to be stochastic and, if $\Gamma$, $\mathcal{G}$,
$\alpha_t$ and $\mathbb{Q}$ are properly built, it is given by the original Markov chain.
In this case, we say that 
$(\Gamma,\mathcal{G}, \alpha_t, \mathbb{Q})$ is a dilation of the Markov evolution in $E$.

Actually, as in the discrete time context, given only the state space $E$ (arbitrary but finite), we build a 
\emph{universal} dilation
$(\Gamma,\mathcal{G},\alpha_t,\{\mathbb{Q}\})$, where $\{\mathbb{Q}\}$ is an entire family of distributions 
which can produce any Markov chain in $E$: 
every Markov chain can be dilated by taking always the same model
$(\Gamma,\mathcal{G},\alpha_t)$ for the environment and the global evolution, and by choosing every time the proper distribution $\mathbb{Q}$
for the initial state of the environment.
Moreover, not only our construction allows to interpret each Markov chain as the stochastic
dynamics resulting from the coupling with an environment, but at the same time it also represents the chain
via an innovation process, 
where the innovation now is dynamically provided by the environment.

Our aim is similar to the aim of Lewis and Maassen \cite{LM84} when they consider classical
mechanics and, taken a linear Hamiltonian system modelling a particle and its
environment, they describe how Gibbs states of the whole system lead to stationary Gaussian stochastic
processes for the observables pertaining to the particle under consideration. However, we do not look
for good global states, but for good states $\mathbb{Q}$ of the environment alone which lead to Markov
evolutions of the system $E$, our particle, for every independent choice of its initial state.

The second aim of the paper is to show that such dilations in Classical
Probability are really analogous to the quantum dilations which inspire them. We prove that every CMS in $E$, considered on any finite time interval, admits an extension to a QDS for 
which we can find a quantum dilation which is itself an extension of the classical dilation $(\Gamma,\mathcal{G},\alpha_t,\mathbb{Q})$ that we build for the CMS.

Thus these dilations are explicit constructions of that classical
structures which would appear by studying the abelian algebras left invariant by quantum stochastic flows
(\textsection 4.27 in \cite{ALV}).

However, we shall not embed a whole universal dilation
$(\Gamma,\mathcal{G},\alpha_t,\{\mathbb{Q}\})$ in the quantum world, as quantum dilations do not exhibit the same
universality and they strictly depend on the QDS under consideration, so that it is not enough
to change the environment state to get another QDS.

In the sequel, given a complex function $f$ on a domain $E$, we shall denote with the same symbol $f$
also its extension on a domain $E\times\Gamma$, $f(i,\gamma)=f(i)$. Similarly, given a map $\phi:E\to E$, we
shall denote with the same symbol $\phi$ also its extension, by tensorizing with the identity, on a domain
$E\times\Gamma$ to $E\times\Gamma$, $\phi(i,\gamma)=\big(\phi(i),\gamma\big)$.

\section{Preliminaries}
We consider a system with finite state space $E=\{1,\ldots,N\}$ and power $\sigma$-algebra $\mathcal{E}$, 
fixed for the whole paper.
We denote by $P=(P_{ij})_{i,j\in E}$ a stochastic matrix in $E$, so that $P_{ij}\geq0$ and $\sum_j
P_{ij}=1$ for every $i$. As usual, we identify the elements of the complex abelian $*$-algebra 
$\mathcal{L}^\infty(\mathcal{E})$, the system random variables
$f:E\to\mathbb{C}$, with the column vectors in $\mathbb{C}^N$, so that every $P$ 
defines an operator in $\mathcal{L}^\infty(\mathcal{E})$, which
describes the one-step evolution $f\mapsto Pf$, $f\in\mathcal{L}^\infty(\mathcal{E})$.
We denote by $D$ a deterministic matrix in $E$, that is a stochastic matrix with
a 1 in each row. Every $D$ describes with matrix terminology a deterministic evolution $\beta$, where
\begin{equation}\label{detmatrix}
    D=(D_{ij})_{i,j\in E}, \qquad \qquad \beta:E\to E, \qquad \qquad D_{ij}=\delta_{\beta(i),j},
\end{equation}
so that $Df=f\circ\beta$.
The invertible maps in $E$ correspond to the special cases of permutation matrices. Labelled all
deterministic matrices and the corresponding maps with indexes $\ell$ belonging to $L=\{1,\ldots,N^N\}$,
every stochastic matrix $P$ determines the weights
\begin{equation}\label{qdcc}
p_\ell=P_{1\beta_\ell(1)}\cdots P_{N\beta_\ell(N)}, \qquad \ell\in L,
\end{equation}
which give a probability on the power $\sigma$-algebra of $L$ and provide the representation
\begin{equation}\label{dcc}
    P=\sum_{\ell\in L}p_\ell \, D_\ell, \qquad \qquad 
    p_\ell\geq0, \quad \sum_{\ell\in L} p_\ell=1.
\end{equation}

We denote by $R=(R_{ij})_{i,j\in E}$ a transition rate matrix in $E$, so that $R_{ij}\geq0$ for every 
$i\neq j$ and $\sum_j R_{ij}=0$ for every $i$. Every $R$ generates a CMS $\e^{Rt}$, $t\geq0$, which consists
of stochastic matrices. It describes the continuous-time homogeneous evolution
\begin{equation}\label{cobsev0}
    f\mapsto \e^{Rt}f, \qquad f\in \mathcal{L}^\infty(\mathcal{E}),\;t\geq0.
\end{equation}
For every CMS $\e^{Rt}$, there exists a Markov chain
$\big(\Omega,\mathcal{F},(\mathcal{F}_t)_{t\geq0},(X_t)_{t\geq0},(\mathbb{P}_k)_{k\in E}\big)$ 
with transition probability functions given by $\e^{Rt}$, i.e. a continuous-time
stochastic process of random variables
$X_t:\Omega\to E$, adapted to a filtration $(\mathcal{F}_t)_{t\geq0}$, and a family of probability
measures $\mathbb{P}_k$, $k\in E$, such that the starting distribution of the process depends on $k$,
$X_0$ has Dirac distribution $\delta_k$ under $\mathbb{P}_k$, but the process always enjoys the 
Markov property with transition matrices $\e^{Rt}$:
\begin{equation*}
    \mathbb{P}_k(X_{t+s}=j|\mathcal{F}_t) = \mathbb{P}_k(X_{t+s}=j|X_t)
    = (\e^{Rs})_{X_tj}, \qquad \forall k,j\in E,\;t,s\geq0.
\end{equation*}
Thus, a system random variable $f\in \mathcal{L}^\infty(\mathcal{E})$ has now a stochastic evolution
given by the $*$-unital homomorphism
\begin{equation}\label{stsystev}
    j_t:\mathcal{L}^\infty(\mathcal{E})\to\mathcal{L}^\infty(\mathcal{F}_t), \qquad f\mapsto
    j_t(f):=f(X_t), \qquad t\geq0,
\end{equation}
and the evolution \eqref{cobsev0} admits the representation
\begin{equation}\label{cobsev1}
    \Big(\e^{Rt}f\Big) (k) = \mathbb{E}_k\big[f(X_t)\big], 
    \qquad \forall f\in \mathcal{L}^\infty(\mathcal{E}).
\end{equation}

By means of the uniformization technique, a continuous-time Markov chain can always be realized as a
discrete-time Markov chain moved by an independent Poisson process. Given the transition rate matrix $R$,
taken a rate $\lambda>0$ and a stochastic matrix $P$ such that
\begin{equation}\label{RPl}
R=\lambda(P-\bbbone),
\end{equation}
taken a 
discrete-time Markov chain $X^*_n$ with transition matrix $P$ and an
independent Poisson process $N(t)$ with rate $\lambda$, then the process $X_t=X^*_{N(t)}$ is a 
continuous-time Markov chain with transition rate matrix $R$ (e.g.\ \cite{EK}).
For example, one can take
\begin{equation}\label{lP}
\lambda=\max_{i\in E}\{-R_{ii}\}, \qquad\qquad P= \bbbone + \frac{1}{\lambda}\,R.
\end{equation}

Moreover, by means of representation \eqref{dcc} for $P$, 
it is always possible to realize the Markov chain via an
innovation process. Let $\Gamma_{\mathbb{R}_+}^L$ be the canonical space of
a marked simple point process on $\mathbb{R}_+$ with mark space $L$, that is the set of sequences
$\gamma=(\ell_n,t_n)_{n\in\mathbb{N}}$ where $\ell_n\in L$, $t_n\in[0,+\infty]$, $0< t_1\leq
t_2\leq\ldots\leq+\infty$, and $t_n<t_{n+1}$ for all $t_n<+\infty$ \cite{BB}. 
Denoting by $\Gamma_{\mathbb{R}_+}$ the canonical space of a simple point process on $\mathbb{R}_+$, we have
$\Gamma_{\mathbb{R}_+}^L=L^{\mathbb{N}}\times\Gamma_{\mathbb{R}_+}$.
Then one can take
\begin{gather}
    \nonumber \Omega=E\times L^{\mathbb{N}}\times\Gamma_{\mathbb{R}_+},\qquad 
    \omega=\big(i,(\ell_n)_{n\in\mathbb{N}},(t_n)_{n\in\mathbb{N}}\big),\\
    \nonumber X_0(\omega)=X^*_0(\omega)=i,\qquad Y_n(\omega)=\ell_n,\qquad
    T_n(\omega)=t_n,\qquad X^*_n=\beta_{Y_n}(X^*_{n-1}),\qquad n\in\mathbb{N} \\
    \label{Dmc} N_\ell(t)=\sum_{n\in\mathbb{N}} I_{(Y_n=\ell,\;T_n\leq t)},\qquad
    N(t)=\sum_{\ell\in L}N_\ell(t),\qquad X_t=X^*_{N(t)}, \qquad t\geq0, \\
    \nonumber \mathcal{F}=\sigma(X_0)\otimes\sigma(Y_n;\;n\in\mathbb{N})\otimes\sigma(T_n;\;n\in\mathbb{N}),
    \qquad \mathcal{F}_t=\sigma\big(X_0, N_\ell(s);\;\ell\in L,\;0\leq s\leq t \big),\\
    \nonumber \mathbb{P}_k=\delta_{k}\otimes p^{\otimes\mathbb{N}}\otimes Q_{\mathbb{R}_+}^\lambda,
\end{gather}
where we denote by $I_A$ the indicator of an event $A$, by $p$ the probability on $L$ associated to $P$, 
and by $Q_{\mathbb{R}_+}^\lambda$ the probability on
$\sigma(T_n;\;n\in\mathbb{N})$ such that the random variables $T_n$ are the arrival times of a Poisson process
with rate $\lambda$. Then the random variables $Y_n$ are i.i.d.\ with distribution 
$p$, the counting
processes $N_\ell(t)$ are independent Poisson processes with rates $p_\ell\lambda$, the counting process
$N(t)$ is a Poisson process with rate $\lambda$, the process $X^*_n$ is a discrete-time Markov chain with
transition matrix $P$ and $X_t$ is a continuous-time Markov chain with transition rate matrix $R$. The random
variables $Y_n$ are called marks. For every $t\geq0$, the chain state $X_t$ is a deterministic function of
$X_0$ and of $N_\ell(s)$, $\ell\in L$, $0\leq s\leq t$, and the random process
$\big(N_\ell(t)\big)_{\ell\in L}$ is an innovation process for $X_t$.

\section{Dilations and universal dilations}

\paragraph{Dilation of a classical Markov semigroup.}
We call dilation of the CMS $\e^{Rt}$ in $\mathcal{L}^\infty(\mathcal{E})$ a term
\begin{equation*}
    \Big(\Omega,\mathcal{F},(\mathcal{F}_t)_{t\geq0},(Z_t)_{t\geq0},
    (\mathbb{P}_{k})_{k\in E}\Big)
\end{equation*}
where
\begin{itemize}
\item every $Z_t=(X_t,\Upsilon_t)$ is a random variable on $(\Omega, \mathcal{F})$ with values in 
$(E\times\Gamma,\mathcal{E}\otimes\mathcal{G})$, being
$(\Gamma,\mathcal{G})$ a fixed measurable space,
\item the term
$\big(\Omega, \mathcal{F}, (\mathcal{F}_t)_{t\geq0}, (X_t)_{t\geq0}, (\mathbb{P}_{k})_{k\in
E}\big)$ is a Markov chain with transition rates $R$,
\item the random variable $(X_0,\Upsilon_0)$ has distribution $\delta_k\otimes\mathbb{Q}$ under 
$\mathbb{P}_{k}$,
being $\mathbb{Q}$ a fixed distribution on $\mathcal{G}$,
\item there exists one-parameter group $(\alpha_t)_{t\in\mathbb{R}}$ of measurable maps 
$\alpha_t:E\times\Gamma\to E\times\Gamma$ such that $\alpha_0=\operatorname{Id}$ and 
$Z_t=\alpha_t(Z_0)$ for every $t\geq0$.
\end{itemize}

Thus, besides the system $E$, a second system is introduced, an environment with state space
$(\Gamma,\mathcal{G})$. Their states $X_t$ and $\Upsilon_t$ are asked to be random variables on a 
same measurable
space $(\Omega,\mathcal{F})$ such that the global state $Z_t=(X_t,\Upsilon_t)$ undergoes a 
deterministic invertible homogeneous evolution $\alpha_t$. Therefore all the $X_t$ and $\Upsilon_t$ 
are determined by $Z_0$, so that $X_t$ and $\Upsilon_t$ are measurable with respect to
$\sigma(Z_0)=\sigma(X_0,\Upsilon_0)\subseteq\mathcal{F}$ and, depending on the probability 
chosen on $\mathcal{F}$, they are deterministic if and only if $Z_0$ is. Nevertheless, a probability 
$\mathbb{P}_{k}$ typically fixes only the value of $X_0$.
The space $(\Omega, \mathcal{F})$ is also endowed with a filtration $\mathcal{F}_t$. 
Note that only the $X_t$ are asked to be adapted to $\mathcal{F}_t$ so that, in particular,
$\Upsilon_0$ does not have to be $\mathcal{F}_0$-measurable. Therefore
the $X_t$ are not trivially $\mathcal{F}_0$-measurable, even if their values 
are completely determined by the values of $X_0$ and 
$\Upsilon_0$, and,
neglecting the environment, each $\big(\Omega, \mathcal{F}, (\mathcal{F}_t)_{t\geq0}, (X_t)_{t\geq0}, 
\mathbb{P}_{k}\big)$ can be a non trivial stochastic process. What we ask
is that $\big(\Omega, \mathcal{F}, (\mathcal{F}_t)_{t\geq0}, (X_t)_{t\geq0}, 
\mathbb{P}_{k}\big)$ actually is a Markov chain starting from $k$ with transition rates 
$R$. At the same time however, this Markov chain is compatible 
with a deterministic, invertible and homogeneous model for the evolution of $E$ coupled with an
environment $\Gamma$. In particular, as $X_0=k$, the whole stochasticity of the process is due only 
to the randomness of the unobserved initial state $\Upsilon_0$ of the environment. 

A dilation gives another interpretation of every evolution \eqref{cobsev0}, 
compatible with \eqref{cobsev1}:
\begin{equation*}
    \Big(e^{Rt}f\Big)(k) = \mathbb{E}_{k}\big[f(X_t)\big] 
    = \mathbb{E}_{k}\big[f(Z_t)\big] = 
    \mathbb{E}_{k}\Big[f\big(\alpha_t(k,\Upsilon_0)\big)\Big], \qquad \forall f\in
    \mathcal{L}^\infty(\mathcal{E}).
\end{equation*}
Indeed, the stochastic evolution \eqref{stsystev} of a system variable $f\in \mathcal{L}^\infty(\mathcal{E})$ 
is now described by the $*$-unital homomorphism
\begin{equation}\label{stsystev2}
    j_t:\mathcal{L}^\infty(\mathcal{E})\to\mathcal{L}^\infty(\mathcal{F}_t), \qquad     
    f\mapsto j_t(f):=f(X_t)=f(Z_t)=f\circ\alpha_t(Z_0),
\end{equation}
which is injective as $\alpha_t$ is invertible. And now we could also consider global
random variables $F:E\times\Gamma\to\mathbb{C}$ and their evolution
$F\mapsto F(Z_t)=F\circ\alpha_t(Z_0)$.

\paragraph{Universal dilation.}
Let us denote by $\mathcal{R}$ the set of transition rate matrices $R$ in $E$.
We call universal dilation of the CMS's in $\mathcal{L}^\infty(\mathcal{E})$ a term
\begin{equation*}
    \Big(\Omega,\mathcal{F},(\mathcal{F}_t)_{t\geq0},(Z_t)_{t\geq0},
    (\mathbb{P}_{k,R})_{k\in E,R\in\mathcal{R}}\Big)
\end{equation*}
such that every $\big(\Omega,\mathcal{F},(\mathcal{F}_t)_{t\geq0},(Z_t)_{t\geq0},
(\mathbb{P}_{k,R})_{k\in E}\big)$ is a dilation of the corresponding semigroup $\e^{Rt}$.
We call universal such a dilation because we ask that the same 
$\Omega$, $\mathcal{F}$, $\mathcal{F}_t$ and $Z_t$ allow
to represent all the CMS's in $\mathcal{L}^\infty(\mathcal{E})$, 
with the change of the probabilities
$\mathbb{P}_{k,R}$ alone. Therefore, both the environment state space $(\Gamma,\mathcal{G})$ and 
the global evolution
$\alpha_t$ depend only on the state space $E$, not on the particular CMS to be dilated.

\paragraph{Poisson dilation and Poisson universal dilation.}
In order to show that every state space $E$ admits a universal dilation, now we consider a particular
classes of dilations and of universal dilations. Let us describe all the special requirements we are
interested in for the dilation of a CMS $\e^{Rt}$. 

First of all we want the sample space $\Omega$ to be just $E\times\Gamma$, the state space of the 
global system. As we want it to describe all the possible initial global states, we ask
the random variable $Z_0$ to be the identity function and $X_0$ and $Y_0$ to be the coordinate 
variables: if $\omega=(i,\gamma)$, then $Z_0(\omega)=\omega$, $X_0(\omega)=i$ and
$\Upsilon_0(\omega)=\gamma$. Thus, for all $t\geq0$, $Z_t=(X_t,\Upsilon_t)=\alpha_t\circ
Z_0=Z_0\circ\alpha_t$, $X_t=X_0\circ\alpha_t$ and $\Upsilon_t=\Upsilon_0\circ\alpha_t$. 

We are interested in an
environment with state space $\Gamma$ equal to $\ess$, the canonical space of a marked simple point process on
$\mathbb{R}$ with finite mark space $G$, that is the set of sequences
$\gamma=(g_n,t_n)_{n\in\mathbb{Z}}$ where $g_n\in G$, $t_n\in[-\infty,+\infty]$, $-\infty\leq\ldots\leq
t_{-1}\leq t_0\leq0< t_1\leq t_2\leq\ldots\leq+\infty$, and $t_n<t_{n+1}$ for all $|t_n|<\infty$.
Denoting by $\Gamma_\mathbb{R}$ the canonical space of a simple point process on $\mathbb{R}$, we have
$\Gamma_\mathbb{R}^G=G^{\mathbb{Z}}\times\Gamma_\mathbb{R}$. Moreover,
with clear meaning of symbols, $\ess=\essn\times\essp$.
Thus the environment state $\gamma$ is a whole trajectory of a marked simple point process.
Later, the introduction of the global evolution $\alpha_t$ will allow the following rough interpretation of
$\gamma$ and of its time-parameter. If the environment state at time 0 is $\gamma=(g_n,t_n)_{n\in\mathbb{Z}}$,
then at every instant $t_n>0$ it will provide a sudden shock of type $g_n$ to the system, thus causing an
instantaneous transition which will be determined by the system state and by $g_n$.

We explicitly introduce the marks $Y_n(\omega)=g_n$, the arrival
times $T_n(\omega)=t_n$, the processes $N_g(t)=\sum_{n\in\mathbb{N}}I_{(Y_n=g,\;0<T_n\leq t)}$, counting the arrivals of type $g$,
and the process $N(t)=\sum_{g\in G}N_g(t)$, counting all the arrivals.
Then, endowed $G$ with its power $\sigma$-algebra, we want all the functions so far introduced to be
measurable and so we ask $\mathcal{G}$ to be the natural $\sigma$-algebra 
$\mathcal{G}_{\mathbb{R}}:=\sigma(\Upsilon_0)=\sigma(Y_n;\;n\in\mathbb{Z})\otimes\sigma(T_n;\;n\in\mathbb{Z})$ 
on $\ess=G^{\mathbb{Z}}\times\Gamma_\mathbb{R}$, and $\mathcal{F}$ 
to be $\mathcal{E}\otimes\mathcal{G}_{\mathbb{R}}=\sigma(X_0,\Upsilon_0)$ on $\Omega$. 
Supposing that at time 0
only $X_0$ is observed and that later only the information carried by the processes
$N_g(t)$ is acquired, we want the filtration $\mathcal{F}_t=\sigma(X_0, N_g(s);\; g\in G,\; 0\leq s\leq t)$.

In order to get consistence between these definitions and the global evolution $\alpha_t$, we ask that
\begin{equation}\label{ge}
\alpha_t=\begin{cases} \vartheta_t\circ\psi_t,&\text{if } t\geq0,\\ 
\psi_{|t|}^{-1}\circ\vartheta_{t},&\text{if } t\leq0,\end{cases}
\end{equation}
where $\vartheta_t$ and $\psi_t$ are as follows. The family of maps $\vartheta_t$ is the group of the left 
shifts
\begin{equation}\label{eve}
\vartheta_t:\ess\to\ess,\qquad
\vartheta_t\big((g_n,t_n)_{n\in\mathbb{Z}}\big)=(g_n,t_n-t)_{n\in\mathbb{Z}}, \qquad t\in\mathbb{R},
\end{equation}
where a renumbering is understood if $t_1-t\leq0$. Every $\vartheta_t$ is
extended on $\Omega$ by tensorizing with the identity. The family of maps 
$\psi_t$ has to be a right cocycle w.r.t.\
$\vartheta_t$ giving a global evolution which, up to $\vartheta_t$, simply couples the system with every mark
provided by the environment, always via a same invertible interaction. More precisely, first we ask an
invertible map
\begin{equation}\label{ipc1}
\phi:E\times G\to E\times G,
\end{equation}
which gives the instantaneous coupling between the system and a single mark. 
Thus the coupling between the system and $m$
subsequent marks is given by the invertible map
\begin{equation}\label{ipc2}\begin{split}
\varphi_m&:E\times\Big(\bigtimes_{n=1}^m G_n\Big)\to E\times\Big(\bigtimes_{n=1}^m G_n\Big), \qquad G_n\equiv
G,\\
\varphi_m&=\phi_m\circ\cdots\circ\phi_1, \qquad \phi_n=\phi:E\times G_n\to E\times G_n.
\end{split}\end{equation}
Then we ask that, given $\phi$ and $\varphi_m$, for every $t\geq0$ the coupling between the system and the 
marks provided in the time-interval $(0,t]$ is given by
\begin{equation}\label{ipc3}
\psi_t:E\times\essp\to E\times\essp, \qquad \psi_t\big(i,(g_n,t_n)_{n\in\mathbb{N}}\big)=
\big(i',(g'_n,t'_n)_{n\in\mathbb{N}}\big),
\end{equation}
where
\begin{equation}\label{ipc4}
t_m\leq t<t_{m+1} \quad\Rightarrow\quad \begin{cases}
\big(i',(g'_n)_{n=1}^m\big) = \varphi_m\big(i,(g_n)_{n=1}^m\big),\quad &\\ 
g'_n=g_n,&\text{if } n>m,\\
t'_n=t_n,&\text{if } n\in\mathbb{N}.\end{cases}
\end{equation} 
Extended also $\psi_t$ on $\Omega$,
for every invertible $\phi$ definitions \eqref{ipc1}-\eqref{ipc4} automatically give a family
of bimeasurable maps satisfying the cocycle property
\begin{equation*}
\psi_{t+s}=\vartheta_{-t}\circ\psi_s\circ\vartheta_t\circ\psi_t, 
\qquad \forall t,s\geq0,
\end{equation*}
and so definition \eqref{ge} automatically give a group of measurable maps in $E\times\ess$. 
Moreover, denoting by $\varphi^E_m$ the projection of $\varphi_m$ on $E$, the state of the system
at a positive time $t$ is $X_t=X_0\circ\alpha_t=X_0\circ\psi_t=\varphi_{N(t)}^E(X_0,Y_1,\ldots,Y_{N(t)})$,
which is automatically 
adapted to $\mathcal{F}_t$. Roughly speaking, when the system and the environment evolve from time 0 to time 
$t>0$, first the map $\psi_t$ couples the initial state $X_0$ of the system with the marks 
$Y_1,\ldots,Y_{N(t)}$, giving the
system state $X_t=\varphi_{N(t)}^E(X_0,Y_1,\ldots,Y_{N(t)})$, and then
the shift $\vartheta_t$ prepares the subsequent marks 
%$Y_{N(t)+1}$, $Y_{N(t)+2},\ldots$
for the future couplings with the system. Thus $\vartheta_t$ could be interpreted as a
free evolution of the environment. 

At long last, we consider the probabilities $\mathbb{P}_{k}$. Of course, they have to be factorized as
$\delta_k\otimes\mathbb{Q}$ on
$\mathcal{F}=\mathcal{E}\otimes\mathcal{G}_{\mathbb{R}}$. We also require the distribution
$\mathbb{Q}$ to be factorized as
$q^{\otimes\mathbb{Z}}\otimes Q_\mathbb{R}^\lambda$ on
$\mathcal{G}_{\mathbb{R}}=\sigma(Y_n;\;n\in\mathbb{Z})\otimes\sigma(T_n;\;n\in\mathbb{Z})$, where $q$ is a
probability on $G$, and $Q_\mathbb{R}^\lambda$ is the probability on $\sigma(T_n;\;n\in\mathbb{Z})$
such that the random variables $T_n$ are the arrival times of a Poisson process with some rate $\lambda$. 
Then the marks $Y_n$ are i.i.d.\ with distribution $q$, the counting processes
$N_g(t)$ are independent Poisson processes with rates $q_g\lambda$, and $N(t)$ is a Poisson process with rate
$\lambda$. As we shall verify in the following Proposition, 
this guarantees the Markov property for the process $X_t$ with respect to $\mathcal{F}_t$, so that the only
point is to check if the resulting rates are the desired ones.

A dilation like this will be called Poisson in the following.
Summarizing, a dilation is Poisson if 
\begin{itemize}
\item $(\Gamma,\mathcal{G})=(\ess,\mathcal{G}_{\mathbb{R}})$, the canonical space of a marked simple point process on $\mathbb{R}$ with finite
mark space $G$,
\item $\Omega=E\times\ess, \quad \omega=(i,\gamma)=
\big(i,(g_n,t_n)_{n\in\mathbb{Z}}\big)\in\Omega, 
\;\;\; i\in E, \;\;\; \gamma\in\ess, \;\;\; g_n\in G, \;\;\; t_n\in\mathbb{R}$,
\item $X_0(\omega)=i, \quad \Upsilon_0=(Y_n,T_n)_{n\in\mathbb{Z}}, \quad \Upsilon_0(\omega)=\gamma, \quad 
Y_n(\omega)=g_n, \quad T_n(\omega)=t_n, \quad Z_0(\omega)=\omega$,
\item $N_g(t)=\sum_{n\in\mathbb{N}}I_{(Y_n=g,\;0<T_n\leq t)}, \qquad N(t)=\sum_{g\in G}N_g(t)$,
\item $\mathcal{F}=\mathcal{E}\otimes\mathcal{G}_{\mathbb{R}}=
\sigma(X_0)\otimes\sigma(Y_n;\;n\in\mathbb{Z})\otimes\sigma(T_n;\;n\in\mathbb{Z})$, 
\item $\mathcal{F}_t=\sigma(X_0, N_g(s);\;g\in G,\;0\leq s\leq t), \quad t\geq0$,
\item $\mathbb{P}_{k}=\delta_k\otimes \mathbb{Q}=
\delta_k\otimes q^{\otimes\mathbb{Z}}\otimes Q_\mathbb{R}^\lambda$,
\item $\alpha_t$ given by \eqref{ge}, with the shift $\vartheta_t$ given by \eqref{eve} and with a cocycle 
$\psi_t$ given by \eqref{ipc1}-\eqref{ipc4}.
\end{itemize} 
Then, for every $t\geq0$ we have
\begin{itemize}
\item $Z_t=\alpha_t\circ Z_0=Z_0\circ\alpha_t=(X_t,\Upsilon_t)$, 
\item $\Upsilon_t=\Upsilon_0\circ\alpha_t$, 
\item $X_t=X_0\circ\alpha_t=X_0\circ\psi_t=\varphi_{N(t)}^E(X_0,Y_1,\ldots,Y_{N(t)})$.
\end{itemize}
A Poisson dilation is therefore specified by the term $\big(G,\phi,\mathbb{Q}\big)$.

With a Poisson dilation the evolution of every global random variable is described by the group of $*$-automorphisms
\begin{equation*} J_t:\mathcal{L}^\infty(\mathcal{F})\to\mathcal{L}^\infty(\mathcal{F}), \qquad J_t(F) = F\circ\alpha_t, \qquad t\in\mathbb{R} . 
\end{equation*}

Note that a Poisson dilation realizes a Markov chain via an innovation process, with the
innovation dynamically provided by the environment.

Let $\phi^E$ denote the projection of $\phi$ on $E$.

\begin{proposition}\label{PdMp}
Let $\big(\Omega,\mathcal{F},(\mathcal{F}_t)_{t\geq0},(X_t)_{t\geq0},(\mathbb{P}_k)_{k\in E}\big)$ be a
stochastic process given by a Poisson dilation 
$\big(G,\phi,q^{\otimes\mathbb{Z}}\otimes Q_\mathbb{R}^\lambda\big)$
for an arbitrary choice of $G$, $\phi$, $q$ and $\lambda$. Then the process is a Markov chain with transition
rate matrix
\begin{equation}\label{RPl2}
R=\lambda(P-\bbbone),\qquad\qquad P_{ij}=\sum_{g\in G}q_g\,\delta_{\phi^E(i,g),j}.
\end{equation}
\end{proposition}

\noindent {\sl Proof.}
Firstly, we set
\begin{equation*}
X^*_0=X_0,\qquad X^*_n=\phi^E(X^*_{n-1},Y_n), \;\; n\in\mathbb{N}, \qquad
\mathcal{F}^*_m=\sigma(X^*_0,Y_n;\;1\leq n\leq m),\;\; m=0,1,\ldots
\end{equation*}
Then $X_t=X^*_{N(t)}$ and 
$\big(\Omega, \mathcal{F}, (\mathcal{F}^*_n)_{n\geq0}, (X^*_n)_{n\geq0}, (\mathbb{P}_{k})_{k\in
E}\big)$ is a discrete-time Markov chain with transition matrix $P$. Moreover,
for every $k\in E$, the independence of $(Y_n)_{n\in\mathbb{N}}$ and $(T_n)_{n\in\mathbb{N}}$ under 
$\mathbb{P}_{k}$ implies the independence of the processes $X^*_n$ and $N(t)$.

Secondly, by the latter independence and the Markov property of $X^*_n$, we show that for $f\in\mathcal{L}^\infty(\mathcal{E})$, $t\geq0$ and $m=0,1,\ldots$
\begin{equation}\label{passint}
\mathbb{E}_{k}\Big[f(X^*_{N(t)+m})\Big|\mathcal{F}_t\Big]=\big(P^mf\big)(X^*_{N(t)}).
\end{equation}
To see this, let $A\in\sigma\big(N(s);\;0\leq s\leq t\big)$, $B\in\mathcal{F}^*_n$, $n\in\mathbb{N}$. Then
\begin{multline*}
\int_{A\cap B\cap(N(t)=n)}f(X^*_{N(t)+m})\,\dd\mathbb{P}_{k}
=\int_{A\cap B\cap(N(t)=n)}f(X^*_{n+m})\,\dd\mathbb{P}_{k}\\
{}=\mathbb{P}_{k}\Big(A,N(t)=n\Big)\int_B f(X^*_{n+m})\,\dd\mathbb{P}_{k}
=\mathbb{P}_{k}\Big(A,N(t)=n\Big)\int_B \big(P^mf\big)(X^*_n)\,\dd\mathbb{P}_{k}\\
=\int_{A\cap B\cap(N(t)=n)}\big(P^mf\big)(X^*_n)\,\dd\mathbb{P}_{k}.
\end{multline*}
Since $\big\{A\cap B\cap(N(t)=n);\; A\in\sigma\big(N(s);\;0\leq s\leq t\big),\, B\in\mathcal{F}^*_n,\,
n\in\mathbb{N}\big\}$ is closed under finite intersections and generates $\mathcal{F}_t$, Eq.\eqref{passint}
follows by the Dynkin class theorem. Finally, since the Poisson process $N(t)$ has independent increments and
rate $\lambda$, 
\begin{multline*}
\mathbb{E}_{k}\Big[f(X_{t+s})\Big|\mathcal{F}_t\Big]
=\mathbb{E}_{k}\Big[f(X^*_{N(t)+N(t+s)-N(t)})\Big|\mathcal{F}_t\Big]\\
=\sum_{m=0}^\infty\mathbb{E}_{k}\Big[f(X^*_{N(t)+m})\Big|\mathcal{F}_t\Big]\,\mathbb{P}_{k}\Big(N(t+s)-N(t)=m\Big)
=\sum_{m=0}^\infty \big(P^mf\big)(X^*_{N(t)})\,\e^{-\lambda s}\frac{(\lambda s)^m}{m!}\\
=\big(\e^{\lambda(P-\bbbone)s}f\big)(X_t)
\end{multline*}
for all $t,s\geq0$. Hence 
$\big(\Omega, \mathcal{F}, (\mathcal{F}_t)_{t\geq0}, (X_t)_{t\geq0}, (\mathbb{P}_{k})_{k\in
E}\big)$ is a Markov chain in $E$ with transition rate matrix $R=\lambda(P-\bbbone)$.
\QED

A universal dilation will be called Poisson if it is given by a family of Poisson dilations 
$\big(G,\phi,\mathbb{Q}_R\big)$, $R\in\mathcal{R}$, all of them with the same $G$ and $\phi$. Thus we fix the
model for the environment and for the global evolution with $G$ and $\phi$, and then we require the existence 
of a family of initial distributions $\mathbb{Q}_R$ for the environment state, 
each one giving rise to a different CMS for the system. A Poisson universal dilation is therefore specified by
the term $\big(G,\phi,(\mathbb{Q}_R)_{R\in\mathcal{R}}\big)$. 
A Poisson universal dilation is not uniquely determined by the state space $E$, but it always exists.

\begin{theorem}\label{cud}
For every finite state space $E$, there exists a Poisson universal dilation\\
$\big(G,\phi,(\mathbb{Q}_R)_{R\in\mathcal{R}}\big)$ of the classical Markov 
semigroups in $\mathcal{L}^\infty(\mathcal{E})$.
\end{theorem}

\noindent {\sl Proof.}
We only have to exhibit a proper mark space $G$, together with the coupling $\phi$ and the probability 
measures $\mathbb{Q}_R$ on $\mathcal{G}_\mathbb{R}$.

We can take the same space $G$ and coupling $\phi$ used in \cite{bm7} for the analogous result in
discrete-time.
Given $E=\{1,\ldots,N\}$ and the set $L$ labelling the all possible maps $\beta:E\to E$, we set
\begin{equation*}
    G=E\times L, \qquad\qquad
    i,j,k\in E, \qquad \ell\in L, \qquad g=(j,\ell)\in G.
\end{equation*}

Arbitrarily fixed $j=1$, we focus on points $(1,\ell)$ in $G$. Thus, taken two points
$\big(i,(1,\ell)\big)\neq\big(i',(1,\ell')\big)$ in $E\times G$, we get
$\big(\beta_\ell(i),(i,\ell)\big)\neq\big(\beta_{\ell'}(i'),(i',\ell')\big)$ 
and so we can find an invertible map
\begin{equation}\label{coupling}
    \phi:E\times G\to E\times G, \qquad\qquad 
    \phi\big(i,(j,\ell)\big) = \begin{cases}\big(\beta_\ell(i),(i,\ell)\big), \quad & \text{if }j=1,\\ 
    \quad \ldots \;, & \text{if }j\neq1.\end{cases}
\end{equation}
We choose an arbitrary $\phi$ satisfying \eqref{coupling}.

Given a transition rate matrix $R$, we consider the rate $\lambda$ and the stochastic matrix
$P$ associated to $R$ by Eq.\eqref{lP}. Via $P$ and Eq.\eqref{qdcc}, we obtain a probability $p$ on $L$,
that we use to define on $G=E\times L$ the probability
\begin{equation*}
q=\delta_1\otimes p.
\end{equation*}
Thus $\displaystyle P_{ij}=\sum_{\ell\in L}p_\ell\,\delta_{\beta_\ell(i),j}=\sum_{g\in G}q_g\,\delta_{\phi^E(i,g),j}$
and, if we define $\mathbb{Q}_R=q^{\otimes\mathbb{Z}}\otimes Q_\mathbb{R}^\lambda$, then every stochastic process
$\big(\Omega, \mathcal{F}, (\mathcal{F}_t)_{t\geq0}, (X_t)_{t\geq0}, (\mathbb{P}_{k,R})_{k\in E}\big)$ is a 
Markov chain with rates $R$, independently of the definition of $\phi\big(i,(j,\ell)\big)$ for $j\neq1$.
Therefore $\big(G,\phi,(\mathbb{Q}_{R})_{R\in\mathcal{R}}\big)$
is a Poisson universal dilation of the classical Markov semigroups in
$\mathcal{L}^\infty(\mathcal{E})$.
\QED

Similarly to the corresponding construction in discrete-time, even if each $g\in G$
has two components, $g=(j,\ell)$, the probability constructed in the proof is always
concentrated only  on those $g$ of the kind
$g=(1,\ell)$, but we need the first component $j$ to define an
invertible $\phi$. Analogously, we are considering the evolution only for positive
times so that all the $(g_n,t_n)$, $n\leq0$, are never involved in the interaction with the
system, but they are needed to define an invertible shift $\vartheta_t$.

We could find also other probabilities on $\mathcal{F}$, different from $\mathbb{P}_{k,R}$, 
but inducing the same 
law for the process $X_t$. Indeed, not only $(Y_n,T_n)_{n\leq0}$ does not effect the evolution of 
$X_t=X_0\circ\psi_t$, but the representation \eqref{RPl} usually holds for other $\lambda$ and $P$ different
from \eqref{lP}, just as
the representation \eqref{dcc} usually holds also for other probabilities $p$ different from \eqref{qdcc}, and
all these different choices could be as well employed in the construction.

\paragraph{The cocycle approach to Poisson dilations.}
Given a CMS $\e^{Rt}$ with a Poisson dilation 
$\big(G,\phi,q^{\otimes\mathbb{Z}}\otimes Q_\mathbb{R}^\lambda\big)$, where $G$, $\phi$, $\lambda$, $q$ 
can be defined as in the proof of Theorem \ref{cud} or not, we can rewrite the $*$-unital injective 
homomorphism \eqref{stsystev2} as
\begin{equation}\label{clfl}
    j_t:\mathcal{L}^\infty(\mathcal{E})\to\mathcal{L}^\infty(\mathcal{F}_t), \qquad f\mapsto
    j_t(f):=f(X_t)=f\circ\psi_t, \qquad t\geq0.
\end{equation}
If we denote by $\mathbb{E}_g[f\circ\phi]$ the system random variable in $\mathcal{L}^\infty(\mathcal{E})$
defined by $i\mapsto f\circ\phi(i,g)$, then the stochastic evolution \eqref{clfl} satisfies
\begin{equation}\label{clfleq}
    j_0(f)=f(X_0), \qquad j_t(f)=\sum_{g\in G}
    j_{t^-}\Big(\mathbb{E}_g[f\circ\phi]-f\Big)\,\dd N_g(t), \qquad
    \forall f\in\mathcal{L}^\infty(\mathcal{E}), \; t\geq0.
\end{equation}
For every $f\in\mathcal{L}^\infty(\mathcal{E})$, this is a stochastic differential equation for the
$\mathcal{F}_t$-adapted process $j_t(f)$
with respect to the noises $N_g(t)$.

More generally, as long as we consider only system random variables neglecting the 
environment, we can define a Markov chain with rates $R$ avoiding $\alpha_t$ and $\vartheta_t$ and 
employing only the deterministic, invertible, but inhomogeneous global evolution 
$\psi_t$, which never involves the environment
components $(g_n,t_n)$ for $n\leq0$. That means to reduce the sample space $\Omega$
from $E\times\ess$ to $E\times\essp$, to restrict here $\mathcal{F}$,
$\mathcal{F}_t$ and $\mathbb{P}_{k}$, and to define the cocycle $\psi_t$, either by
\eqref{ipc1}-\eqref{ipc4} or by \eqref{ipc1} and \eqref{clfleq}. 
Anyway, thanks to the cocycle properties of $\psi_t$, it is always possible to introduce later $\essn$ and 
the shift $\vartheta_t$, in order to recover the whole environment
state space $(\ess,\mathcal{G}_{\mathbb{R}})$, the evolution $\alpha_t$ and the initial environment
distribution $q^{\otimes\mathbb{Z}}\otimes Q_\mathbb{R}^\lambda$, so that the two constructions are equivalent
and can be considered different descriptions of the same 
situation. 

Choosing the cocycle approach, a Poisson dilation of a CMS $\e^{Rt}$ gives a Markov chain 
$\big(\Omega, \mathcal{F}, (\mathcal{F}_t)_{t\geq0}, (X_t)_{t\geq0}, (\mathbb{P}_{k})_{k\in E}\big)$,
where
\begin{gather}
    \Omega=E\times\essp,\qquad 
    \omega=\big(i,(g_n,t_n)_{n\in\mathbb{N}}\big),\nonumber\\
    X_0(\omega)=X^*_0(\omega)=i,\qquad Y_n(\omega)=g_n,\qquad
    T_n(\omega)=t_n,\qquad X^*_n=\phi^E(X^*_{n-1},Y_n),\qquad n\in\mathbb{N} \nonumber\\
    N_g(t)=\sum_{n\in\mathbb{N}} I_{(Y_n=g,\;T_n\leq t)},\qquad
    N(t)=\sum_{g\in G}N_g(t),\qquad X_t=X^*_{N(t)}, \qquad t\geq0, \label{stPMc}\\
    \mathcal{F}=\sigma(X_0)\otimes\sigma(Y_n;\;n\in\mathbb{N})\otimes\sigma(T_n;\;n\in\mathbb{N}),
    \qquad \mathcal{F}_t=\sigma\big(X_0, N_g(s);\;g\in G,\;0\leq s\leq t \big),\nonumber\\
    \mathbb{P}_{k}=\delta_{k}\otimes q^{\otimes\mathbb{N}}\otimes Q_{\mathbb{R}_+}^\lambda, \nonumber
\end{gather}
and where, again, $R$ is related to $\phi$, $q$ and $\lambda$ by Eq.~\eqref{RPl2}. This chain
is specified by the term 
$\big(G,\phi,q^{\otimes\mathbb{N}}\otimes Q_{\mathbb{R}_+}^\lambda\big)$.

Let us underline that the Markov chain \eqref{stPMc} is similar to \eqref{Dmc}, 
as also this one is represented via a discrete-time Markov chain and an innovation process. 
Nevertheless, this chain is endowed with a reacher structure because the cocycle $\psi_t$
implicitly introduces also the
deterministic invertible homogeneous evolution $\alpha_t$.

\section{Dilations of classical Markov semigroups and of quantum dynamical semigroups}

We want to compare a Poisson dilation with the dilation of a QDS in
Quantum Probability.

Given a Hilbert space $\is$, always complex separable in the paper, 
let us denote its vectors by $h$, or $|h\rangle$ using
Dirac's notation, so that $\langle h'|h\rangle$ denotes the scalar product (linear in $h$) and
$|h'\rangle\langle h|$ denotes the operator $h'' \mapsto \langle h|h''\rangle \,h'$. Let
$\mathcal{B}(\is)$ be the complex $*$-algebra of bounded operators in $\is$.
Let $\Gamma[\is]$ denote the symmetric Fock space over $\is$ and, for every
$h\in\is$, let $e(h)\in\Gamma[\is]$ denote the corresponding normalized exponential vector. Given a 
Hilbert space $L^2(\mu)$, with a probability measure $\mu$ on some measurable space, and given a measurable
complex function $f$ on the same measurable space,
let $m_f$ denote the multiplication operator
\begin{equation*}
m_f:\operatorname{Dom}(m_f)\to L^2(\mu), \qquad \operatorname{Dom}(m_f)=\{h\in L^2(\mu)\;:\; fh\in L^2(\mu)\},
\qquad m_f\,h=fh,
\end{equation*}
which is bounded if and only if $f\in L^\infty(\mu)$.
Given two Hilbert spaces $\is$ and $\mathcal{K}$ and a vector $\kappa\in\mathcal{K}$, let
$\mathbb{E}_\kappa:\mathcal{B}(\is)\otimes\mathcal{B}(\mathcal{K})\to\mathcal{B}(\is)$ denote the conditional
expectation with respect to $|\kappa\rangle\langle\kappa|$.
In the sequel, given an operator $a$ in $\is$, we shall identify it with its extension
$a\otimes\bbbone_\mathcal{K}$ in $\is\otimes\mathcal{K}$.

\paragraph{Quantum extension of a CMS.}
We denote by $\mathcal{L}$ a Lindblad operator in $\mathcal{B}(\is)$, that is an operator
$\mathcal{L}:\mathcal{B}(\is)\to\mathcal{B}(\is)$ admitting the representation
\begin{equation}\label{lindbladian}
\mathcal{L}a=\im [H,a] + \sum_{z\in Z}\Big(R_z^*aR_z - \frac{1}{2}\{R_z^*R_z,a\}\Big), \qquad
    a\in\mathcal{B}(\is),
\end{equation}
where $Z$ is a discrete index set, $H$ and $R_z$ belong to $\mathcal{B}(\is)$, $H^*=H$, $\sum_{z}R_z^*R_z$ 
strongly converges in $\mathcal{B}(\is)$, and where $[\cdot,\cdot]$ and $\{\cdot,\cdot\}$ denote
the commutator and the anticommutator respectively. 
Let us recall that representation \eqref{lindbladian} is not unique,
and that every $\mathcal{L}$ generates a uniformly continuous QDS
$T_t=\e^{\mathcal{L}t}$, $t\geq0$, which consists of bounded, completely positive, normal and identity
preserving operators
$T_t:\mathcal{B}(\is)\to\mathcal{B}(\is)$. Actually Eq.\eqref{lindbladian} gives the complete characterization
of the generator of a uniformly continuous QDS.

In order to extend a CMS in $\mathcal{L}^\infty(\mathcal{E})$ by a QDS in some $\mathcal{B}(\is)$, we take
$\is=L^2(\mu_E)$, where $\mu_E$ is the uniform probability on $(E,\mathcal{E})$, and, denoted by 
$\{|i\rangle\}_{i\in E}$ its canonical basis, we embed $\mathcal{L}^\infty(\mathcal{E})=L^\infty(\mu_E)$ in 
$\mathcal{B}(\is)$ by the $*$-isomorphism $f \mapsto m_f=\sum_{i\in E} f(i) \, |i\rangle\langle i|$
between $L^\infty(\mu_E)$ and the subalgebra of the multiplication operators in $L^2(\is)$.

Since $\is$ is finite dimensional, every QDS in $\mathcal{B}(\is)$ is uniformly continuous.
We say that a QDS $\e^{\mathcal{L}t}$ in $\mathcal{B}(\is)$
extends a CMS $\e^{Rt}$ in $L^\infty(\mu_E)$ if 
\begin{equation*}
    \e^{\mathcal{L}t}\,m_f = m_{\e^{Rt}f}, \qquad \forall f \in L^\infty(\mu_E),\; t\geq0.
\end{equation*}
It is enough to check that $\mathcal{L}m_f = m_{Rf}$ for all $f\in L^\infty(\mu_E)$.
Such extension always exists, but it is not unique at all. For example,
given the transition rate matrix $R$, considered a rate $\lambda$ and a stochastic matrix $P$ associated 
to $R$ via \eqref{RPl}, taken a probability $p$ associated to $P$ via \eqref{dcc}, 
using notation \eqref{detmatrix}, 
the CMS $\e^{Rt}$ is extended by the QDS $\e^{\mathcal{L}t}$ generated by 
\begin{equation}\label{Rext}
    \mathcal{L}a = \lambda \Big( \sum_{\begin{subarray}{c} \scriptscriptstyle\ell\in L \\
    \scriptscriptstyle i\in E \end{subarray}}
    p_\ell\,|i\rangle\langle\beta_\ell(i)|\,a\,|\beta_\ell(i)\rangle\langle i| - a \Big), \qquad
    a\in\mathcal{B}(\is),
\end{equation}
which admits representation \eqref{lindbladian} with $H=0$, $Z=G=E\times L$,
$R_z=R_{(i,\ell)}=\sqrt{p_\ell\lambda}\,|\beta_\ell(i)\rangle\langle i|$.

\paragraph{HP-dilation of a QDS.}
Given a uniformly continuous QDS $\e^{\mathcal{L}t}$ in $\mathcal{B}(\is)$, a typical Quantum Probability 
construction employs a Hudson-Parthasarathy equation to dilate
$\e^{\mathcal{L}t}$ at the same time by a quantum stochastic flow and by a group of $*$-automorphisms
\cite{ALV,F1,F2,HP,M1,M2,M,MS,P}. 

Chosen a representation \eqref{lindbladian} for $\mathcal{L}$, taken the Hilbert space $\ms$ generated by the
basis $\{|z\rangle\}_{z\in Z}$,
let $\dd \Lambda_{zz'}(t)$, $\dd A_z(t)$ and $\dd A^\dagger_z(t)$ be the corresponding canonical quantum 
noises in $\Gamma[L^2(\mathbb{R}_+;\ms)]$.
Fix an arbitrary vector $\nu=\sum_z\nu_z|z\rangle\in\ms$ and an arbitrary unitary operator
$S=\sum_{zz'}S_{zz'}\otimes|z\rangle\langle z'|$
in $\is\otimes\ms$. Define the bounded operators in $\is$
\begin{gather}\nonumber
L_z=R_z-\sum_{z'}S_{zz'}\nu_{z'},\\
\label{HPcdaL}
H_0=H+\frac{\im}{2}\sum_{zz'}\Big(R_z^*S_{zz'}\nu_{z'}-\bar{\nu}_{z'}S^*_{zz'}R_{z}\Big)=H_0^*.
\end{gather}
Then the quantum stochastic differential equation for adapted processes in 
$\is\otimes\Gamma[L^2(\mathbb{R}_+;\ms)]$
\begin{equation*}
\begin{split}
\dd V_t&=\bigg[\sum_{zz'}\Big(S_{zz'}-\delta_{zz'}\Big)\dd\Lambda_{zz'}(t) 
- \sum_{zz'}L^*_zS_{zz'}\dd A_{z'}(t)
+ \sum_zL_z\dd A^\dagger_z(t) - \Big(\im H_0 +\frac{1}{2}\sum_zL^*_zL_z\Big)\dd t\bigg]V_t,\\
V_0&=\bbbone,
\end{split}
\end{equation*}
is a Hudson-Parthasarathy equation. The properties of the coefficients guarantee that it admits a unique
solution $V_t$, which is a strongly continuous unitary cocycle. The HP-dilation of $\e^{\mathcal{L}t}$ 
is identified by the triple $(\ms,\nu,V_t)$ as follows. 

The quantum stochastic flow is 
\begin{equation}\label{qsf}
j_t:\mathcal{B}(\is)\to\mathcal{B}(\is)\otimes\mathcal{B}\big(\Gamma[L^2(\mathbb{R}_+;\ms)]\big), \qquad
j_t(a)=V_t^*\,a\,V_t, \qquad t\geq0,
\end{equation}
which satisfies the quantum stochastic differential equation
\begin{multline*}
\dd j_t(a) = \sum_{zz'}j_t\Big(\sum_{z''}S^*_{z''z}aS_{z''z'}-\delta_{zz'}a\Big)\dd\Lambda_{zz'}(t)
+ \sum_z j_t\Big(\sum_{z'}S^*_{z'z}[a,L_{z'}]\Big)\dd A^\dagger_z(t)\\
+ \sum_z j_t\Big(\sum_{z'}[L^*_{z'},a]S_{z'z}\Big)\dd A_z(t)
+ j_t\Big(\im[H_0,a]+\sum_z\big(L^*_zaL_z - \frac{1}{2}\{L^*_zL_z,a\}\big)\Big)\dd t.
\end{multline*}

Denoted by $\Theta_t$ the second quantization of the left shift in $L^2(\mathbb{R};\ms)$, $v(r)\mapsto
v(r+t)$, embedded $\is\otimes\Gamma[L^2(\mathbb{R}_+;\ms)]$ in $\is\otimes\Gamma[L^2(\mathbb{R};\ms)]$ and
extended here every operator by tensorizing with the identities, consider the
strongly continuous unitary group in $\is\otimes\Gamma[L^2(\mathbb{R};\ms)]$
\begin{equation}\label{qSge}
U_t=\begin{cases}\Theta_tV_t,&\text{if }t\geq0,\\V^*_{|t|}\Theta_t,&\text{if }t\leq0.\end{cases}
\end{equation}
The group of $*$-automorphisms is
\begin{equation}\label{qHge}
J_t:\mathcal{B}(\is)\otimes\mathcal{B}\big(\Gamma[L^2(\mathbb{R};\ms)]\big) \to
\mathcal{B}(\is)\otimes\mathcal{B}\big(\Gamma[L^2(\mathbb{R};\ms)]\big), \qquad
J_t(A)=U_t^*\,A\,U_t, \qquad t\in\mathbb{R}.
\end{equation}
Of course $J_t(a)=j_t(a)$ for every $a\in\mathcal{B}(\is)$ and
$t\geq0$.

Then, for every $T>0$, taken $v\in L^2(\mathbb{R}_+;\ms)$ and $u\in L^2(\mathbb{R};\ms)$ such that
$v(t)=u(t)=\nu$ for every $0\leq t\leq T$, both $\big(j_t,e(v)\big)$ and $\big(J_t,e(u)\big)$ dilate 
$e^{\mathcal{L}t}$ in the time interval $[0,T]$, that is
\begin{equation}\label{HPdil}
e^{\mathcal{L}t}a = \mathbb{E}_{e(v)}\Big[j_t(a)\Big] 
= \mathbb{E}_{e(u)}\Big[J_t(a)\Big] \qquad \forall
a\in\mathcal{B}(\is),\;0\leq t\leq T.
\end{equation}
In particular, if $\nu=0$, then the choices $v=0$ and $u=0$ give the usual dilations which hold
for every $t\geq0$.

Let us call such a construction a HP-dilation $(\ms,\nu,V_t)$. 
This is highly non-unique, as representation
\eqref{lindbladian} and the choices of $S$ and $\nu$ are not.

Thanks to the adaptedness of $V_t$, when $0\leq t\leq T<\infty$, we can consider the first equality of
Eq.~\eqref{HPdil} with $V_t$ adapted process in $\is\otimes\Gamma[L^2((0,T);\ms)]$, with the quantum
stochastic flow
$j_t:\mathcal{B}(\is)\to\mathcal{B}(\is)\otimes\mathcal{B}\big(\Gamma[L^2((0,T);\ms)]\big)$ and with the state
$e(v)\in\Gamma[L^2((0,T);\ms)]$.

\paragraph{Quantum extension of a Poisson dilation.}
Given a CMS $\e^{Rt}$, we can either represent it in Classical Probability with a Poisson dilation 
$\big(G,\phi,q^{\otimes\mathbb{Z}}\otimes Q_\mathbb{R}^\lambda\big)$, either extend it to a QDS
$\e^{\mathcal{L}t}$ and then represent this latter in Quantum Probability with a HP-dilation $(\ms,\nu,V_t)$.
Now we want to prove that, if the quantum extension $\mathcal{L}$ and its dilation $(\ms,\nu,V_t)$ are 
properly chosen, then $(\ms,\nu,V_t)$ is a quantum extension of 
$\big(G,\phi,q^{\otimes\mathbb{Z}}\otimes Q_\mathbb{R}^\lambda\big)$. 
This shows that a Poisson dilation is a
classical analogue of a HP-dilation.

The first step to study this relationship is to embed a Poisson dilation
$\big(G,\phi,q^{\otimes\mathbb{Z}}\otimes Q_\mathbb{R}^\lambda\big)$ in the quantum world. 
We want this embedding at the Hilbert space level, as this is the level where
usual quantum stochastic calculus is defined. Therefore we should introduce first a proper measure on
$(\ess,\mathcal{G}_{\mathbb{R}})$. Taken $\mu_E$ on $E$, the proper measure on
$(\ess,\mathcal{G}_{\mathbb{R}})$ should give a product measure on 
$(E\times\ess,\mathcal{E}\otimes\mathcal{G}_{\mathbb{R}})$ invariant for the deterministic invertible
evolutions $\alpha_t$, $\vartheta_t$, $\psi_t$. Thus the natural choice would be the probability measure
$\mu_G^{\otimes\mathbb{Z}}\otimes Q_{\mathbb{R}}^\lambda$, where
$\mu_G$ denotes the uniform probability on $G$. Nevertheless, it would be singular with respect to the initial
environment distribution $q^{\otimes\mathbb{Z}}\otimes Q_\mathbb{R}^\lambda$, so that this latter could not be
obtained from a state in $L^2(\mu_G^{\otimes\mathbb{Z}}\otimes Q_{\mathbb{R}}^\lambda)$.
Therefore now we fix a finite
time horizon $T>0$ and we focus only for $0\leq t\leq T$ on the Markov chain  
$\big(G,\phi,q^{\otimes\mathbb{N}}\otimes Q_{\mathbb{R}_+}^\lambda\big)$.

Using notations \eqref{stPMc}, we introduce the $\sigma$-algebras on $\essp$
\begin{gather*}
\mathcal{G}_{\mathbb{R}_+} = \sigma(Y_n;\;n\in\mathbb{N})\otimes\sigma(T_n;\;n\in\mathbb{N}),\\
\mathcal{G}_t = \sigma\big(N_g(s);\;g\in G,\;0\leq s\leq t\big)\subseteq\mathcal{G}_{\mathbb{R}_+}.
\end{gather*}
Then we consider the probability measures 
$\mathbb{Q}=\mu_G^{\otimes\mathbb{N}}\otimes Q_{\mathbb{R}_+}^\lambda$ on $\mathcal{G}_{\mathbb{R}_+}$ and
$\mathbb{P}=\mu_E\otimes\mathbb{Q}$ on $\mathcal{F}=\mathcal{E}\otimes\mathcal{G}_{\mathbb{R}_+}$. 
Given $T>0$, we introduce also their
restrictions $\mathbb{Q}_T$ and $\mathbb{P}_T$ to $\mathcal{G}_T$ and $\mathcal{F}_T$. Thus the processes
$N_g(t)$, $0\leq t\leq T$, are i.i.d.\ Poisson processes with rates $\lambda/|G|$. Then we consider the
Hilbert space $L^2(\mathbb{P}_T)=L^2(\mu_E)\otimes L^2(\mathbb{Q}_T)$ and we embed
$L^\infty(\mathbb{P}_T)=L^\infty(\mu_E)\otimes L^\infty(\mathbb{Q}_T)$ in
$\mathcal{B}\big(L^2(\mathbb{P}_T)\big)$ by the $*$-isomorphism $F\mapsto m_F$ between
$L^\infty(\mathbb{P}_T)$ and the bounded multiplication operators in $L^2(\mathbb{P}_T)$.

The cocycle $\psi_t$ generates the family of operators
\begin{equation*}
\Psi_t:L^2(\mathbb{P}_T)\to L^2(\mathbb{P}_T), \qquad \Psi_t\,\xi=\xi\circ\psi_t^{-1}, \qquad 0\leq t\leq T.
\end{equation*}
\begin{proposition}Every $\Psi_t$ is a well defined unitary operator in $L^2(\mathbb{P}_T)$.
\end{proposition}

\noindent {\sl Proof.}
Every $\Psi_t$ is well defined and unitary because, thanks to the invertibility of $\phi$, the maps $\psi_t$
and $\psi_t^{-1}$ preserve the probability measure $\mathbb{P}_T$.

To see this, let $A\in\sigma\big(N(s);\;0\leq s\leq T\big)$, $B=(X_0=i,Y_1=g_1,\ldots,Y_m=g_m)$, 
$0\leq n\leq m$. Then
\begin{multline*}
\psi_t\Big(A\cap\big(N(t)=n\big)\cap\big(N(T)=m\big)\cap B\Big)
= A\cap\big(N(t)=n\big)\cap\big(N(T)=m\big)\cap \psi_t(B)=\\
A\cap\big(N(t)=n\big)\cap\big(N(T)=m\big)\cap \Big((X_0, Y_1,\ldots,Y_n)=\varphi_n(i,g_1,\ldots,g_n),
Y_{n+1}=g_{n+1},\ldots,Y_m=g_m\Big),
\end{multline*}
and so, thanks to the invertibility of $\varphi_n$,
\begin{multline*}
\mathbb{P}_T\Big(\psi_t\Big(A\cap\big(N(t)=n\big)\cap\big(N(T)=m\big)\cap B\Big)\Big) 
%\\
%= \mathbb{P}_T\Big(A,N(t)=n,N(T)=m\Big)
%\, \mathbb{P}_T\Big((X_0, Y_1,\ldots,Y_n)=\varphi_n(i,g_1,\ldots,g_n),Y_{n+1}=g_{n+1},\ldots,Y_m=g_m\Big)\\
= Q_{\mathbb{R}_+}^\lambda\Big(A,N(t)=n,N(T)=m\Big)
\, \frac{1}{|E|}\,\frac{1}{|G|^m}\\
%\mu_E\otimes\mu_G^{\otimes\mathbb{N}}
%\Big((X_0, Y_1,\ldots,Y_n)=\varphi_n(i,g_1,\ldots,g_n),Y_{n+1}=g_{n+1},\ldots,Y_m=g_m\Big)
= \mathbb{P}_T\Big(A,N(t)=n,N(T)=m,B\Big).
\end{multline*}
Since these events form a set closed under finite intersections and generate $\mathcal{F}_t$,
the equality $\mathbb{P}_T\big(\psi_t(A)\big)=\mathbb{P}_T(A)$ follows for all $A\in\mathcal{F}_T$. Analogously, 
$\mathbb{P}_T\big(\psi_t^{-1}(A)\big)=\mathbb{P}_T(A)$ for all $A\in\mathcal{F}_T$.
\QED

\noindent The unitary family $\Psi_t$ is a quantum extension of the deterministic invertible evolution 
$\psi_t$, as
\begin{equation*}
\Psi_t^*\,m_F\,\Psi_t=m_{F\circ\psi_t}, \qquad \forall F\in L^\infty(\mathbb{P}_T), \;0\leq t\leq T.
\end{equation*}

In order to study the quantum stochastic calculus properties of the family of operators $\Psi_t$, we need to
introduce a good isomorphism between $L^2(\mathbb{Q}_T)$ and a Fock space $\Gamma[L^2((0,T);\ms)]$.
We take $\ms=L^2(\mu_G)$
with its canonical basis $\{|g\rangle\}_{g\in G}$ and we choose the isomorphism which
diagonalizes the number processes $\Lambda_{gg}(t)$ and which maps 
the state $e(v_0)$, with $v_0(t)=\sum_{g}\sqrt{\lambda/|G|}\,|g\rangle$ for all $0\leq t\leq T$,
to the constant function $1\in L^2(\mathbb{Q}_T)$. 
Denoted by $\mu_T$ the Lebesgue measure on $(0,T)$, let us
make the identification $L^2\big((0,T);\ms\big)=L^2(\mu_G\otimes\mu_T)$ so that every $v:(0,T)\mapsto\ms$, 
$v(t)=\sum_g v_g(t) |g\rangle$, corresponds to $v:G\times(0,T)\to\mathbb{C}$, $v(g,t)=\sqrt{|G|}\,v_g(t).$
Thus we introduce the isomorphism
\begin{gather*}
I_T:\Gamma[L^2((0,T);\ms)]\to L^2(\mathbb{Q}_T),\\
I_T[e(v)]=\exp\Big(-\frac{1}{2}\|v\|^2+\frac{1}{2}\lambda T\Big)
\prod_{n=1}^{N(T)}\frac{v(Y_n,T_n)}{\sqrt{\lambda}}, \qquad v\in L^2(\mu_G\otimes\mu_T)\cap
\mathcal{C}\big([0,T];\ms\big).
\end{gather*}
The operator $I_T$ is an isometry which turns out to be unitary thanks to the chaotic representation property
of the Poisson process. Then
\begin{equation*}
I_T\,\Lambda_{gg}(t)\,I_T^{-1}=m_{N_g(t)}, \qquad\forall g\in G,\;0\leq t\leq T,
\end{equation*}
and $I_T[e(v_0)]=1$. Any choice of a state
$\eta\in\Gamma[L^2((0,T);\ms)]$ with $\eta\neq e(v_0)$ corresponds to a
change of probability on $\mathcal{G}_T$, that is the choice of a probability with Radon
Nikodym derivative $|I_T[\eta]|^2$ with respect to $\mathbb{Q}_T$.

Extended $I_T$ to an isomorphism between $\is\otimes\Gamma[L^2((0,T);\ms)]$ and $L^2(\mathbb{P}_T)$,
we define the family of unitary operators
\begin{equation}\label{psiext}
V_t:=I_T^{-1}\,\Psi_t\,I_T:\is\otimes\Gamma[L^2((0,T);\ms)]\to\is\otimes\Gamma[L^2((0,T);\ms)], \qquad 
0\leq t\leq T.
\end{equation}
Thus the unitary quantum evolution in 
$\mathcal{B}(\is)\otimes\mathcal{B}\big(\Gamma[L^2((0,T);\ms)]\big)$, $A\mapsto V_t^*\,A\,V_t$, 
$0\leq t\leq T$, admits an invariant abelian subalgebra where it gives just the classical evolution 
$\circ\psi_t$ in $L^\infty(\mathbb{P}_T)$. 
Indeed, the algebra of the multiplication operators 
$\{m_F,\;F\in L^\infty(\mathbb{P}_T)\} \subseteq \mathcal{B}\big(L^2(\mathbb{P}_T)\big)$ is $*$-isomorphic 
to the algebra $\{M_F:=I_T^{-1}\,m_F\,I_T,\;F\in L^\infty(\mathbb{P}_T)\} 
\subseteq \mathcal{B}(\is)\otimes\mathcal{B}\big(\Gamma[L^2((0,T);\ms)]\big)$, where we have
\begin{equation*}
V_t^*\,M_F\,V_t = M_{F\circ\psi_t}, \qquad \forall F\in L^\infty(\mathbb{P}_T), \; 0\leq t\leq T.
\end{equation*}
In particular, the flow $j_t$ associated to $V_t$ via \eqref{qsf} is a quantum 
extension of the classical homomorphism $j_t$ associated to $\psi_t$ via \eqref{clfl}.

\begin{theorem}\label{cqd} Let $\e^{Rt}$ be a classical Markov semigroup in a finite state space $E$ and let \eqref{stPMc}
be the Markov chain provided by a Poisson dilation 
$(G,\phi,q^{\otimes\mathbb{N}}\otimes Q_{\mathbb{R}_+}^\lambda)$. 
Let $\is=L^2(\mu_E)$, let $\ms=L^2(\mu_G)$, and let 
$V_t$ be the family of unitary operators \eqref{psiext} in $\is\otimes\Gamma[L^2((0,T);\ms)]$. Then
\begin{itemize}
\item[(1)] $V_t$ is a strongly continuous adapted process satisfying the Hudson-Parthasarathy equation
\begin{equation}\label{nHPeq}\begin{split}
\dd V_t&=\sum_{g,g'\in G}\Big(S_{gg'}-\delta_{gg'}\Big)V_t\,\dd\Lambda_{gg'}(t),\qquad 0\leq t\leq T,\\
V_0&=\bbbone,
\end{split}\end{equation}
where $S$ is the unitary operator $\displaystyle\sum_{\begin{subarray}{c} \scriptscriptstyle i\in E\\
\scriptscriptstyle g\in G\end{subarray}}|\phi(i,g)\rangle\langle i,g|$ in
$\is\otimes\ms$;
\item[(2)] the quantum stochastic flow $j_t(a)=V_t^*\,a\,V_t$ satisfies, on the algebra of the
multiplication operators, the quantum stochastic differential equation
\begin{equation}\label{nFeqm}
\dd j_t(m_f) = 
\sum_{g\in G}j_t\Big( \sum_{i\in E}|i\rangle\langle\phi^E(i,g)|\,m_f\,|\phi^E(i,g)\rangle\langle i| 
- m_f \Big)\dd\Lambda_{gg}(t),\qquad\forall f\in L^\infty(\mu_E);
\end{equation}
\item[(3)] taken $v(t)=\nu$ for all $0\leq t\leq T$, with $\nu=\sum_{g\in G}\sqrt{\lambda
q_g}\,|g\rangle\in\ms$, 
\begin{itemize}
\item $\displaystyle |I_T[e(v)]|^2=\frac{\dd q^{\otimes\mathbb{N}}\otimes
Q_{\mathbb{R}_+}^\lambda|_{\mathcal{G}_T}}{\dd\mathbb{Q}_T}$,
\item $(j_t,e(v))$ dilates in $0\leq t\leq T$ the quantum dynamical semigroup $\e^{\mathcal{L}t}$, extending $\e^{Rt}$, which is generated by
\begin{equation}\label{Rext0}
    \mathcal{L}a = \lambda \Big( \sum_{g,g',g''\in G}
    \sqrt{q_{g'}\,q_{g''}}\,S_{gg'}^*\,a\,S_{gg''} - a \Big), \qquad
    a\in\mathcal{B}(\is).
\end{equation}
\end{itemize}
\end{itemize}
\end{theorem}

\noindent Before of the proof, let us write explicitly the relation between $\phi$ and coefficients in 
\eqref{nHPeq}, that is
\begin{equation*}
S_{gg'}=\sum_{i,j\in E}|i\rangle\langle i,g|\phi(j,g')\rangle\langle j|.
\end{equation*}

\noindent {\sl Proof.}
(1) The family of unitary operators $V_t$ is an adapted process in $\is\otimes\Gamma[L^2((0,T);\ms)]$ because
every $V_t$ belongs to
$\mathcal{B}\Big(\is\otimes\Gamma[L^2((0,t);\ms)]\Big)\otimes\bbbone_{\Gamma[L^2((t,T);\ms)]}$. Indeed, for
every $h,h'\in L^2(\mu_E)$, $v,v'\in L^2((0,T);\ms)$, $0\leq t\leq T$,
\begin{multline*}
\langle h\otimes e(v)\pmb{\big|}V_t\,h'\otimes e(v')\rangle \hfill \\
= \exp\Bigl(\lambda T-\frac{1}{2}\|v\|^2-\frac{1}{2}\|v'\|^2\Bigr)\hfill \\
\hfill \cdot
\mathbb{E}_T\!\left[\overline{h}(X_0)\left(\prod_{n=1}^{N(T)}\frac{\overline{v}(Y_n,T_n)}{\sqrt{\lambda}}\right)
\left(\left\{h'(X_0)\,\prod_{n=1}^{N(t)}\frac{v'(Y_n,T_n)}{\sqrt{\lambda}}\right\}\circ\psi_t^{-1}\right)
\prod_{n=N(t)+1}^{N(T)}\frac{v'(Y_n,T_n)}{\sqrt{\lambda}}\right]\\
= \exp\Bigl(\lambda T-\frac{1}{2}\|v\|^2-\frac{1}{2}\|v'\|^2\Bigr) \hfill\\
\cdot
\mathbb{E}_T\!\left[\overline{h}(X_0)\left(\prod_{n=1}^{N(t)}\frac{\overline{v}(Y_n,T_n)}{\sqrt{\lambda}}\right)
\left(\left\{h'(X_0)\,\prod_{n=1}^{N(t)}\frac{v'(Y_n,T_n)}{\sqrt{\lambda}}\right\}\circ\psi_t^{-1}\right)\right]
\\
\hfill\cdot
\mathbb{E}_T\!\left[\prod_{n=N(t)+1}^{N(T)}\frac{\overline{v}(Y_n,T_n)\,v'(Y_n,T_n)}{\lambda}\right]\\
= \langle h\otimes e(v|_{(0,t)})\pmb{\big|}V_t\,h'\otimes e(v'|_{(0,t)})\rangle
\;\langle e(v|_{(t,T)})\pmb{\big|}e(v'|_{(t,T)})\rangle, \hfill
\end{multline*}
where the last $V_t$ is the operator in $\is\otimes\Gamma[L^2((0,t);\ms)]$ defined by Eq.~\eqref{psiext} with 
$T=t$, which a posteriori is identified with its extension in $\is\otimes\Gamma[L^2((0,T);\ms)]$ for $T\geq
t$.

The adapted process $V_t$ is strongly continuous. Since $V_t$ is unitary for every $t$, it is enough to prove
weak continuity. To see this, let $\xi\in L^2(\mathbb{P}_T)$, $\xi'\in L^\infty(\mathbb{P}_T)$, $0\leq t\leq
T$, $0\leq t+s\leq T$, $\Delta N=|N(t+s)-N(t)|$. Then
\begin{multline*}
\left|\mathbb{E}_T\!\left[\overline\xi\;\Psi_{t+s}\,\xi'\right] 
- \mathbb{E}_T\!\left[\overline\xi\;\Psi_t\,\xi'\right]\right|
\leq \int_\Omega |\xi|\cdot\left|\xi'\circ\psi_{t+s}^{-1} - \xi'\circ\psi_t^{-1}\right|\,\dd\mathbb{P}_T\\
= \int_{\Delta N\geq1} |\xi|\cdot\left|\xi'\circ\psi_{t+s}^{-1} - \xi'\circ\psi_t^{-1}\right|\,
\dd\mathbb{P}_T
\leq 2\,\|\xi'\|_{\infty}\int_{\Delta N\geq1} |\xi|\,\dd\mathbb{P}_T \xrightarrow[s\to0]{}0,
\end{multline*}
because $\psi_{t+s}^{-1}(\omega)=\psi_t^{-1}(\omega)$ for all $\omega\in(\Delta N=0)$ and $\mathbb{P}_T(\Delta
N\geq1)=1-\e^{-\lambda|s|}\to0$ as $s\to0$. Thus $\mathbb{E}_T[\overline\xi\,\Psi_t\,\xi']$ is continuous for
$0\leq t\leq T$ if $\xi\in L^2(\mathbb{P}_T)$ and $\xi'\in L^\infty(\mathbb{P}_T)$, and the same is true for 
$\xi\in L^\infty(\mathbb{P}_T)$ and $\xi'\in L^2(\mathbb{P}_T)$, as the same argument works for
$\Psi_t^*\,\xi=\xi\circ\psi_t$. If both $\xi$ and $\xi'$ belong to $L^2(\mathbb{P}_T)$, then the continuity of
$\mathbb{E}_T[\overline\xi\,\Psi_t\,\xi']$ follows by standard arguments taking a sequence 
$\xi_n\in L^\infty(\mathbb{P}_T)$ such that $\xi_n\to\xi$ in $L^2(\mathbb{P}_T)$. Therefore $\Psi_t$ and
$V_t$ are strongly continuous.

In order to show that the unitary strongly continuous adapted process $V_t$ satisfies the Hudson-Parthasarathy
equation \eqref{nHPeq}, it is enough to show that
\begin{multline}\label{wnHPeq}
\mathfrak{f}(t):=\langle h\otimes e(v)\pmb{\big|}V_t\,h'\otimes e(v')\rangle\\
{}= \mathfrak{f}(0)
+ \int_0^t\langle h\otimes e(v)\pmb{\big|}\sum_{g,g'}\Big(S_{gg'}-\delta_{gg'}\Big)V_s\,\overline{v}_{g'}(s)\,v'_g(s)\,
h'\otimes e(v')\rangle\,\dd s,
\end{multline}
for all $h,h'\in L^2(\mu_E)$, $v,v'\in L^2((0,T);\ms)$, $0\leq t\leq T$. To prove this, we compute the right
derivative of $\mathfrak{f}(t)$ in the special case $v,v'\in\mathcal{C}([0,T];\ms)$. Let $\xi=I_T[h\otimes e(v)]$,
$\xi'=I_T[h'\otimes e(v')]$, $0\leq t<t+s\leq T$, $\Delta N=N(t+s)-N(t)$, $\Delta N_g=N_g(t+s)-N_g(t)$. Then
\begin{multline*}
\mathfrak{f}(t+s)-\mathfrak{f}(t) \\
= \mathbb{E}_T\!\left[\left(\overline\xi\circ\psi_{t+s} - \overline\xi\circ\psi_{t}\right)\xi'\right]\hfill\\
= \sum_g \mathbb{E}_T\!\left[\left.\left(\overline\xi\circ\psi_{t+s} - \overline\xi\circ\psi_{t}\right)\xi' 
\right| \Delta N=\Delta N_g=1\right] \mathbb{P}_T(\Delta N=\Delta N_g=1) \hfill\\
\hfill + \mathbb{E}_T\!\left[\left.\left(\overline\xi\circ\psi_{t+s} - \overline\xi\circ\psi_{t}\right)\xi' 
\right| \Delta N\geq2\right] \mathbb{P}_T(\Delta N\geq2)\\
= A + B. \hfill
\end{multline*}
Setting
\begin{equation*}
C_1 = \exp\Bigl(\lambda T-\frac{1}{2}\|v\|^2-\frac{1}{2}\|v'\|^2\Bigr) \, \|h\|_\infty \, \|h'\|_\infty,
\qquad
C_2 = \frac{\|v\|_\infty \, \|v'\|_\infty}{\lambda},
\end{equation*}
the second addendum $B$ is bounded as follows:
\begin{multline*}
|B|\leq \mathbb{E}_T\!\left[\left.\left(|\overline\xi\circ\psi_{t+s}| + |\overline\xi\circ\psi_{t}|\right)
|\xi'| \right| \Delta N\geq2\right] \mathbb{P}_T(\Delta N\geq2) \\
\leq 2C_1\,\mathbb{E}_T\!\left[\left.C_2^{N(T)}\right| \Delta N\geq2\right] \mathbb{P}_T(\Delta N\geq2)\\
= 2 C_1\,\e^{(C_2-1)\lambda T}\left(1-\e^{-C_2\lambda s} - C_2\,\lambda s\,\e^{-C_2\lambda s}\right)
= o(s), \quad \text{as }s\to0^+.
\end{multline*}
In order to deal with the first addendum $A$, for $n\geq1$ let us explicitly introduce the random variables
$Y_n^{(t)}=Y_n\circ\psi_t$, the $n$-th marks at time $t$. Denoting by $\varphi_n^{G_n}$ the projection of the
map $\varphi_n$ on $G_n$ (see Eq.~\eqref{ipc2}), we have
\begin{equation*}
Y_n^{(t)}=Y_n\circ\psi_t
=\begin{cases}\varphi^{G_n}_n(X_0,Y_1,\ldots,Y_n),&\text{if }T_n\leq t,\\
Y_n,&\text{if }T_n>t,\end{cases}
\end{equation*}
so that $Y_n^{(t+s)}=Y_n^{(t)}$ if $T_n\leq t<t+s$, and $Y_n^{(t+s)}=Y_n^{(t)}=Y_n$ if $t<t+s<T_n$.
Denoting by $\phi^G$ the projection of the map $\phi$ on $G$, we also have
$Y_n^{(t+s)}=\phi^G(X_t,Y_n)$ if $T_{n-1}\leq t<T_n\leq t+s$. Then, for every $g\in G$,
\begin{multline*}
\mathbb{E}_T\!\left[\left.\left(\overline\xi\circ\psi_{t+s} - \overline\xi\circ\psi_{t}\right)\xi' 
\right| \Delta N=\Delta N_g=1\right] 
\hfill
\\
= \exp\Bigl(\lambda T-\frac{1}{2}\|v\|^2-\frac{1}{2}\|v'\|^2\Bigr) \,
\mathbb{E}_T\!
\left[
\left(
\prod_{n=N(t+s)+1}^{N(T)}\frac{\overline{v}(Y_n,T_n)\,v'(Y_n,T_n)}{\lambda}
\right)
\right.
\hfill
\\
\cdot
\left(
\overline{h}(\phi^E(X_t,g))\,\frac{\overline{v}(\phi^G(X_t,g),T_{N(t+s)})}{\sqrt{\lambda}}
-\overline{h}(X_t)\,\frac{\overline{v}(g,T_{N(t+s)})}{\sqrt{\lambda}}
\right)
h'(X_0)\,\frac{v'(g,T_{N(t+s)})}{\sqrt{\lambda}}
\\
\hfill
\cdot
\left.
\left.
\prod_{n=1}^{N(t)}\frac{\overline{v}(Y_n^{(t)},T_n)\,v'(Y_n,T_n)}{\lambda}
\right|
\Delta N=\Delta N_g=1
\right]
\\
= \exp\Bigl(\lambda T-\frac{1}{2}\|v\|^2-\frac{1}{2}\|v'\|^2\Bigr) \,
\mathbb{E}_T\!\left[\prod_{n=N(t+s)+1}^{N(T)}\frac{\overline{v}(Y_n,T_n)\,v'(Y_n,T_n)}{\lambda}\right]
\hfill
\\
\cdot
\mathbb{E}_T\!
\left[
\left(
\overline{h}(\phi^E(X_t,g))\,\frac{\overline{v}(\phi^G(X_t,g),T_{N(t+s)})}{\sqrt{\lambda}}
-\overline{h}(X_t)\,\frac{\overline{v}(g,T_{N(t+s)})}{\sqrt{\lambda}}
\right)
h'(X_0)\,\frac{v'(g,T_{N(t+s)})}{\sqrt{\lambda}}
\right.
\\
\hfill
\cdot
\left.
\left.
\prod_{n=1}^{N(t)}\frac{\overline{v}(Y_n^{(t)},T_n)\,v'(Y_n,T_n)}{\lambda}
\right|
\Delta N=\Delta N_g=1
\right].
\end{multline*}
Thus
\begin{multline*}
A
= \exp\Bigl(\lambda T-\frac{1}{2}\|v\|^2-\frac{1}{2}\|v'\|^2\Bigr) \,
\mathbb{E}_T\!\left[\prod_{n=N(t+s)+1}^{N(T)}\frac{\overline{v}(Y_n,T_n)\,v'(Y_n,T_n)}{\lambda}\right]
\\
\cdot
\sum_g
\left\{
\mathbb{E}_T\!
\left[
\left(
\overline{h}(\phi^E(X_t,g))\,\frac{\overline{v}(\phi^G(X_t,g),t)}{\sqrt{\lambda}}
-\overline{h}(X_t)\,\frac{\overline{v}(g,t)}{\sqrt{\lambda}}
\right)
h'(X_0)\,\frac{v'(g,t)}{\sqrt{\lambda}}
\right.
\right.
\hfill
\\
\hfill
\cdot
\left.
\left.
\prod_{n=1}^{N(t)}\frac{\overline{v}(Y_n^{(t)},T_n)\,v'(Y_n,T_n)}{\lambda}
\right|
\Delta N=\Delta N_g=1
\right]
\\
+
\mathbb{E}_T\!
\left[
\left(
\frac{\overline{v}(\phi^G(X_t,g),T_{N(t+s)})\,v'(g,T_{N(t+s)})}{\lambda}
- \frac{\overline{v}(\phi^G(X_t,g),t)\,v'(g,t)}{\lambda}
\right)
\overline{h}(\phi^E(X_t,g))\,h'(X_0)
\right.
\hfill
\\
\hfill
\cdot
\left.
\left.
\prod_{n=1}^{N(t)}\frac{\overline{v}(Y_n^{(t)},T_n)\,v'(Y_n,T_n)}{\lambda}
\right|
\Delta N=\Delta N_g=1
\right]
\\
+
\mathbb{E}_T\!
\left[
\left(
\frac{\overline{v}(g,t)\,v'(g,t)}{\lambda}
- \frac{\overline{v}(g,T_{N(t+s)})\,v'(g,T_{N(t+s)})}{\lambda}
\right)
\overline{h}(X_t)\,h'(X_0)
\right.
\hfill
\\
\hfill
\cdot
\left.
\left.
\left.
\prod_{n=1}^{N(t)}\frac{\overline{v}(Y_n^{(t)},T_n)\,v'(Y_n,T_n)}{\lambda}
\right|
\Delta N=\Delta N_g=1
\right]
\right\}
\frac{\lambda s}{|G|}\,\e^{-\lambda s}
\\
= A_1 + A_2 + A_3,
\hfill
\end{multline*}
where the three addenda $A_l$ correspond to the three expectations enclosed in curly brackets.
Thanks to the continuity of $v$ and $v'$, every function $\overline{v}(g,r)\,v'(g',r)/\lambda$ is continuous with respect to $r$, so that there exists a constant $C_3(s)$ such that
\begin{equation*}
\left|\frac{\overline{v}(g,r)\,v'(g',r)}{\lambda} - \frac{\overline{v}(g,t)\,v'(g',t)}{\lambda}\right| \leq C_3(s), \qquad \forall t\leq r\leq t+s,\; g,g'\in G,
\end{equation*}
with $C_3(s)\to0$ as $s\to0^+$. Since $t<T_{N(t+s)}\leq t+s$, both $A_2$ and $A_3$ are bounded as follows:
\begin{multline*}
|A_l| \leq C_1\,\mathbb{E}_T\!\left[C_2^{N(T)-N(t+s)}\right]\,C_3(s)\,\mathbb{E}_T\!\left[C_2^{N(t)}\right]\,\lambda s \,\e^{-\lambda s} \\
= C_1\,\e^{(C_2-1)\lambda(T-s)}\,C_3(s)\,\lambda s\,\e^{-\lambda s}\, = o(s), \qquad \text{as }s\to0^+.
\end{multline*}
Finally, let us show that $A_1/s$ converges to the right limit. The following conditional expectation has an 
$\mathcal{F}_t$-measurable argument and hence it satisfies
\begin{multline*}
\mathbb{E}_T\!
\left[
\left(
\overline{h}(\phi^E(X_t,g))\,\frac{\overline{v}(\phi^G(X_t,g),t)}{\sqrt{\lambda}}
-\overline{h}(X_t)\,\frac{\overline{v}(g,t)}{\sqrt{\lambda}}
\right)
h'(X_0)\,\frac{v'(g,t)}{\sqrt{\lambda}}
\right.
\hfill
\\
\hfill
\cdot
\left.
\left.
\prod_{n=1}^{N(t)}\frac{\overline{v}(Y_n^{(t)},T_n)\,v'(Y_n,T_n)}{\lambda}
\right|
\Delta N=\Delta N_g=1
\right]
\\
{} = \mathbb{E}_T\!
\left[
\left(
\overline{h}(\phi^E(X_t,g))\,\frac{\overline{v}(\phi^G(X_t,g),t)}{\sqrt{\lambda}}
-\overline{h}(X_t)\,\frac{\overline{v}(g,t)}{\sqrt{\lambda}}
\right)
h'(X_0)\,\frac{v'(g,t)}{\sqrt{\lambda}}
\right.
\hfill
\\
\hfill
\cdot
\left.
\prod_{n=1}^{N(t)}\frac{\overline{v}(Y_n^{(t)},T_n)\,v'(Y_n,T_n)}{\lambda}
\right]
\\
{} = \mathbb{E}_T\!
\left[
\left(
\overline{h}(\phi^E(X_0,g))\,\frac{\overline{v}(\phi^G(X_0,g),t)}{\sqrt{\lambda}}
-\overline{h}(X_0)\,\frac{\overline{v}(g,t)}{\sqrt{\lambda}}
\right)
\left(
\prod_{n=1}^{N(t)}\frac{\overline{v}(Y_n,T_n)}{\sqrt\lambda}
\right)
\right.
\hfill
\\
\hfill
\cdot
\left.
\frac{v'(g,t)}{\sqrt{\lambda}}
\left\{
h'(X_0)
\prod_{n=1}^{N(t)}\frac{v'(Y_n,T_n)}{\sqrt\lambda}
\right\}
\circ\psi^{-1}_t
\right]
\\
{} = \mathbb{E}_T\!
\left[
\left(
\sum_{j,g'}\langle j,g'|\phi(X_0,g)\rangle
\overline{h}(j)\,\frac{\overline{v}(g',t)}{\sqrt{\lambda}}
-\overline{h}(X_0)\,\frac{\overline{v}(g,t)}{\sqrt{\lambda}}
\right)
\left(
\prod_{n=1}^{N(t)}\frac{\overline{v}(Y_n,T_n)}{\sqrt\lambda}
\right)
\right.
\hfill
\\
\hfill
\cdot
\left.
\frac{v'(g,t)}{\sqrt{\lambda}}
\left\{
h'(X_0)
\prod_{n=1}^{N(t)}\frac{v'(Y_n,T_n)}{\sqrt\lambda}
\right\}
\circ\psi^{-1}_t
\right]
\\
{}=\frac{1}{\lambda}\,
\langle
\left(
\sum_{g'} S^*_{g'g}\,v(g',t)-v(g,t)
\right)
h\otimes e(v|_{(0,t)})
\pmb{\big|}
v'(g,t)\,V_t\,h'\otimes e(v'|_{(0,t)})
\rangle
\hfill
\\
\hfill
\cdot
\exp\Bigl(-\lambda t+\frac{1}{2}\|v|_{(0,t)}\|^2+\frac{1}{2}\|v'|_{(0,t)}\|^2\Bigr).
\end{multline*}
Then
\begin{multline*}
\frac{A_1}{s} = 
\exp\Bigl(\lambda T-\frac{1}{2}\|v\|^2-\frac{1}{2}\|v'\|^2\Bigr) \,
\mathbb{E}_T\!\left[\prod_{n=N(t+s)+1}^{N(T)}\frac{\overline{v}(Y_n,T_n)\,v'(Y_n,T_n)}{\lambda}\right]
\hfill
\\
\cdot
\sum_g
\mathbb{E}_T\!
\left[
\left(
\overline{h}(\phi^E(X_t,g))\,\frac{\overline{v}(\phi^G(X_t,g),t)}{\sqrt{\lambda}}
-\overline{h}(X_t)\,\frac{\overline{v}(g,t)}{\sqrt{\lambda}}
\right)
h'(X_0)\,\frac{v'(g,t)}{\sqrt{\lambda}}
\right.
\\
\hfill
\cdot
\left.
\left.
\prod_{n=1}^{N(t)}\frac{\overline{v}(Y_n^{(t)},T_n)\,v'(Y_n,T_n)}{\lambda}
\right|
\Delta N=\Delta N_g=1
\right]
\,\frac{\lambda}{|G|}\,\e^{-\lambda s}
\\
{}=\langle
h\otimes e(v|_{(0,t)})
\pmb{\big|}
\sum_{gg'}\left(S_{gg'}-\delta_{gg'}\right)\,\overline{v}_g(t)\,v'_{g'}(t)\,V_t\,h'\otimes e(v'|_{(0,t)})
\rangle
\,
\langle e(v|_{(t+s,T)})\pmb{\big|}e(v'|_{(t+s,T)})\rangle
\hfill
\\
\hfill
\cdot
\exp\Bigl(-\frac{1}{2}\|v|_{(t,t+s)}\|^2-\frac{1}{2}\|v'|_{(t,t+s)}\|^2\Bigr)
\\
\xrightarrow[s\to0^+]{}
\langle
h\otimes e(v)
\pmb{\big|}
\sum_{gg'}\left(S_{gg'}-\delta_{gg'}\right)\,\overline{v}_g(t)\,v'_{g'}(t)\,V_t\,h'\otimes e(v')
\rangle.
\hfill
\end{multline*}
Therefore, for every $h,h'\in\is$ and every $v,v'\in\mathcal{C}([0,T],\ms)$, the function $\mathfrak{f}$ is a continuous function with continuous right derivative $D_+\mathfrak{f}(t) = \langle h\otimes e(v)|\sum_{g,g'}\left(S_{gg'}-\delta_{gg'}\right)\,\overline{v}_g(t)\,v'_{g'}(t)\,V_t\,h'\otimes e(v')\rangle$. Then $\mathfrak{f}$ is continuously differentiable with $\dd\mathfrak{f}/\dd t=D_+\mathfrak{f}$. Finally, the density of $\mathcal{C}([0,T],\ms)$ in $L^2([0,T],\ms)$ gives Eq.~\eqref{wnHPeq} on the whole exponential domain.

(2) The Hudson-Parthasarathy equation \eqref{nHPeq} for the process $V_t$ determines the quantum stochastic differential equation $\dd j_t(a) = \sum_{g,g'}j_t\Big(\sum_{g''}S^*_{g''g}aS_{g''g'}-\delta_{gg'}a\Big)\dd\Lambda_{gg'}(t)$ for the quantum stochastic flow $j_t(a)=V_t^*\,a\,V_t$. If $a=m_f$ with $f\in L^\infty(\mu_E)$, then
\begin{multline*}
\sum_{g''}S^*_{g''g}\,m_f\,S_{g''g'} = \sum_{g''}\sum_{i,j,i',j'} |i\rangle\langle \phi(i,g)|j,g''\rangle\langle j| \,m_f\, |i'\rangle\langle i',g''|\phi(j',g')\rangle\langle j'|\\
= \delta_{gg'}\sum_{i} |i\rangle\langle \phi^E(i,g)|\,m_f\,|\phi^E(i,g)\rangle\langle i|,
\end{multline*}
so that Eq.˜\eqref{nFeqm} follows for every multiplication operator in $L^2(\mu_E)$.

(3) If we take $\nu=\sum_{g}\sqrt{\lambda q_g}\,|g\rangle\in\ms$ and $v(t)=\nu$ for every $0\leq t\leq T$, then
\begin{equation*}
|I_T[e(v)]|^2=\prod_{n=1}^{N(T)}|G|\,q_{Y_n} = \frac{\dd q^{\otimes\mathbb{N}}\otimes
Q_\lambda|_{\mathcal{G}_T}}{\dd\mathbb{Q}_T}.
\end{equation*}
Moreover, since $V_t$ satisfies a Hudson-Parthasarathy equation with only $\dd\Lambda_{gg'}$ terms, the couple
$(j_t,e(v))$ defines a QDS $\e^{\mathcal{L}t}$ in $\mathcal{B}(\is)$ with Lindblad generator $\mathcal{L}$
admitting representation \eqref{lindbladian} with coefficients $R_g$ and $H$ given by \eqref{HPcdaL} in the 
case $L_g=0$ and $H_0=0$. Thus
\begin{equation*}R_g=\sum_{g'}S_{gg'}\nu_{g'},\qquad\qquad H=0,
\end{equation*}
which give, together with our choice of $\nu$, the Lindblad operator \eqref{Rext0}.
Then
\begin{equation*}
    \mathcal{L}m_f = \lambda \Big( \sum_{\begin{subarray}{c} \scriptscriptstyle g\in G \\
    \scriptscriptstyle i\in E \end{subarray}}
    q_g\,|i\rangle\langle\phi^E(i,g)|\,m_f\,|\phi^E(i,g)\rangle\langle i| - m_f \Big) = m_{Rf}, \qquad \forall f\in L^\infty(\mu_E).
\end{equation*}
\QED

Theorem \ref{cqd} allows to state that the cocycle $\psi_t$ of a Poisson dilation
is a classical analogue of a Hudson-Parthasarathy cocycle, as it admits a quantum extension $V_t$ satisfying
the Hudson-Parthasaraty equation \eqref{nHPeq}. In particular the associated quantum stochastic flow $j_t$ 
satisfies, on the abelian algebra of the multiplication operators in $L^2(\mathbb{P}_T)$, 
the quantum stochastic differential equation \eqref{nFeqm}, which is just a reformulation in operator
terminology of the stochastic differential equation \eqref{clfleq} satisfied by the corresponding classical 
homomorphism $j_t$.

Moreover, Theorem \ref{cqd} allows to state that the whole Poisson dilation 
$(G,\phi,q^{\otimes\mathbb{N}}\otimes Q_{\mathbb{R}_+}^\lambda)$ of a CMS $\e^{Rt}$ is a classical analogue of a HP-dilation $(\ms,\nu,V_t)$ of a QDS $\e^{\mathcal{L}t}$, as there is also a quantum environment state $e(v)$ which gives on side the right initial distribution of the classical environment to dilate $\e^{Rt}$, and on the other side the right QDS to extend $\e^{Rt}$.

Let us note also that if $(G,\phi,q^{\otimes\mathbb{N}}\otimes Q_{\mathbb{R}_+}^\lambda)$ is built as in Theorem \ref{cud}, then
\begin{equation*}
S_{(i,\ell)(1,\ell')}=\delta_{\ell\ell'}\,|\beta_\ell(i)\rangle\langle i|,
\end{equation*}
and the Lindblad operator \eqref{Rext0} becomes just the Lindblad operator \eqref{Rext}.

All of these results are obtained for an arbitrary but finite time horizon $T>0$. They are consistent 
with respect to $T$, but it is not possible to set $T=+\infty$, as there is no isomorphism $I_{\infty}$,
and $v(t)=\nu$ for all $t>0$ does not belong to $L^2(\mathbb{R}_+;\ms)$. Thus, starting from $\psi_t$ and
$V_t$, we can introduce separately the groups $\vartheta_t$ and $\alpha_t$ and the groups $\Theta_t$ and
$U_t$, but we do not have a Hilbert space isomorphism to show that the group of 
$*$-automorphisms \eqref{qHge} gives the group of $*$-automorphisms $\circ\alpha_t$ on an invariant abelian 
subalgebra of multiplication operators.

In order to avoid the finite time horizon and to find a correspondence between the two dilations holding for 
all times, one can look at a Poisson dilation under other isomorphisms different from $I_T$. 
For example, one can consider the usual isomorphism
\begin{gather*}
\widehat{I}:\Gamma[L^2(\mathbb{R}_+;\ms)]\to L^2(\mathbb{Q}),\\
\widehat{I}[e(v)]=\exp\left(-\frac{1}{2}\|v\|^2-\sqrt{\frac{\lambda}{|G|}}\,\sum_{g\in G}\int_{\mathbb{R}_+}
v_g(s)\,\dd s\right) \prod_{n=1}^{\infty}\left(1+\frac{v(Y_n,T_n)}{\sqrt{\lambda}}\right), \quad v\in
\mathcal{C}_\mathrm{c}\big([0,+\infty);\ms\big), 
\end{gather*}
such that
\begin{equation*}
\widehat{I}^{-1}\,m_{N_g(t)}\,\widehat{I}=\Lambda_{gg}(t)+\sqrt{\frac{\lambda}{|G|}}\,\Big(A_g^\dagger(t)+A_g(t)\Big)+\frac{\lambda}{|G|}\,t,
\qquad\forall g\in G,\;0\leq t<\infty, 
\end{equation*}
and $\widehat{I}[e(0)]=1$. This choice leads now to a strongly continuous adapted process of unitary operators
$\widehat{V}_t=\widehat{I}^{-1}\,\Psi_t\,\widehat{I}$ in $\is\otimes\Gamma[L^2(\mathbb{R}_+;\ms)]$ defined for
all $t>0$. Nevertheless the new Hudson-Parthasarathy equation is not as simple as \eqref{nHPeq} and it is
anyhow related to \eqref{nHPeq} by a Weyl transformation. Moreover, still no 
state $\eta\in\Gamma[L^2(\mathbb{R}_+;\ms)]$ can give 
$\displaystyle |\widehat{I}[\eta]|^2=\dd q^{\otimes\mathbb{N}}\otimes Q_{\mathbb{R}_+}^\lambda/\dd\mathbb{Q}$
on $\mathcal{G}_{\mathbb{R}_+}$, and again we can recover the right semigroups $\e^{Rt}$ and
$\e^{\mathcal{L}t}$ only up to a finite time $T$, namely by choosing an environment state $e(v)$,
\begin{equation*}
v(t)=\begin{cases}\sum_{g\in G}\Big(\sqrt{q_g\lambda}-\sqrt{\frac{\lambda}{|G|}}\Big)\,|g\rangle, & \text{if
}0\leq t\leq T,\\ u(t),& \text{if }t>T,\end{cases} 
\end{equation*}
with an arbitrary $u\in L^2((T,+\infty);\ms)$.

In order to eliminate the finite time horizon one can also leave the Hilbert space approach and study this
correspondence on some $C^*$-algebras of bounded operators, but then the connection with quantum stochastic
calculus is less direct. Anyway the basic result remains Theorem \ref{cqd}, which can be employed to find
the preferred isomorphism.  

Let us conclude by remarking that, if Theorem \ref{cqd} allows to interpret a Poisson dilation as a classical
analogue of a HP-dilation, at the same time Eq.~\eqref{nHPeq} is only a particular case of a Hudson-Parthasarathy
equation, so that it also suggests that other classical analogues could be found by coupling the system $E$ 
with processes different from counting processes.

\end{document}